\documentclass[%
 aip,
 amsmath,amssymb,
 reprint,%
]{revtex4-1}

\usepackage{graphicx}% Include figure files
\usepackage{dcolumn}% Align table columns on decimal point
\usepackage{bm}% bold math
\usepackage[utf8]{inputenc}
\usepackage[T1]{fontenc}
\usepackage{mathptmx}
\usepackage{etoolbox}
\usepackage{subcaption}
\usepackage[dvipsnames]{xcolor}
\usepackage{amsthm}
\usepackage{lipsum}
\usepackage{overpic}
\usepackage{algorithm,algcompatible}
\usepackage[hidelinks, pdftex,
            pdfauthor={Rudy Geelen, Laura Balzano, Stephen Wright, Karen Willcox},
            pdftitle={Learning physics-based reduced-order models from data using nonlinear manifolds},
            pdfsubject={Data-driven model reduction, Nonlinear manifolds, Operator inference, Proper orthogonal decomposition},
            pdfkeywords={Data-driven model reduction, Nonlinear manifolds, Operator inference, Proper orthogonal decomposition}]{hyperref}

\definecolor{pblue}{RGB}{8,57,110}
\definecolor{pblue2}{RGB}{135,190,217}
\definecolor{pwhite}{RGB}{248,248,248}
\definecolor{pred2}{RGB}{238,152,120}
\definecolor{pred}{RGB}{106,1,20}
\DeclareMathAlphabet{\mathcal}{OMS}{cmsy}{m}{n}

\usepackage{tikz}
\usepackage{pgfplots}
\pgfplotsset{compat=1.15, legend image code/.code={
\draw[mark repeat=2,mark phase=2]
plot coordinates {
(0cm,0cm)
(0.15cm,0cm)        %% default is (0.3cm,0cm)
(0.3cm,0cm)         %% default is (0.6cm,0cm)
};%
}}
\usetikzlibrary{calc, colorbrewer, arrows}
\usepgfplotslibrary{fillbetween}

\theoremstyle{definition}
\newtheorem{exmp}{Example}[section]

\DeclareMathOperator*{\argmin}{arg\,min}
\DeclareMathOperator{\sech}{sech}

\makeatletter
\def\@email#1#2{%
 \endgroup
 \patchcmd{\titleblock@produce}
  {\frontmatter@RRAPformat}
  {\frontmatter@RRAPformat{\produce@RRAP{*#1\href{mailto:#2}{#2}}}\frontmatter@RRAPformat}
  {}{}
}%
\makeatother

\begin{document}

%\preprint{AIP/123-QED}

\title{Learning physics-based reduced-order models from data using nonlinear manifolds}

\email{rudy.geelen@austin.utexas.edu}

\author{Rudy Geelen}
\affiliation{Oden Institute for Computational Engineering and Sciences, University of Texas at Austin, Austin, TX 78712}%

\author{Laura Balzano}%
\affiliation{Electrical Engineering and Computer Science Department, University of Michigan, Ann Arbor, MI 48109}%

\author{Stephen Wright}%
\affiliation{Computer Sciences Department, University of Wisconsin, Madison, WI 53706}%

\author{Karen Willcox}
\affiliation{Oden Institute for Computational Engineering and Sciences, University of Texas at Austin, Austin, TX 78712}%

\date{\today}% It is always \today, today,

\begin{abstract}
We present a novel method for learning reduced-order models of dynamical systems using nonlinear manifolds. First, we learn the manifold by identifying nonlinear structure in the data through a general representation learning problem. The proposed approach is driven by embeddings of low-order polynomial form. A projection onto the nonlinear manifold reveals the algebraic structure of the reduced-space system that governs the problem of interest. The matrix operators of the reduced-order model are then inferred from the data using operator inference. Numerical experiments on a number of nonlinear problems demonstrate the generalizability of the methodology and the increase in accuracy that can be obtained over reduced-order modeling methods that employ a linear subspace approximation.
\end{abstract}

\maketitle

\begin{quotation}
Model reduction rests on the fundamental assumption that system states in complex, physics-based models can often be represented with a smaller number of variables without a significant loss of information. The identification of such intrinsic, low-dimensional structure in these problems, and the subsequent inference of projection-based reduced-order models lies at the heart of this paper. We treat the construction of nonlinear state approximations of polynomial form as a general representation learning problem. By leveraging data-driven operator inference, we can then learn reduced-order models directly from available snapshot data. These models are physics-informed in that their algebraic structure is dictated by the original, high-dimensional problem. The proposed nonlinear model reduction method is interpretable and effective for reducing large-scale, dynamical-system models.
\end{quotation}

\section{Introduction and background}
\label{sec:introduction}

Projection-based model reduction is comprised of a family of methods that build approximations of complex physics-based models with (generally speaking) orders-of-magnitude reduction in computational complexity. 
Through identification of inherent low-dimensional structure, the cost of many computational tasks can be lowered significantly. 
Model reduction makes tractable many applications in control, uncertainty quantification, optimal experimental design, and inverse problems.\cite{ghattas_willcox_2021} The key idea behind many model reduction methods is to identify a low-dimensional representation for a set of training snapshots by applying data compression. 
The computation of high-dimensional states is then replaced with identification of the coefficients of a basis expansion in the reduced subspace. 
As the effectiveness of the reduction step hinges on the ability to find a sufficiently accurate reduced-dimensional representation of full-state vectors, the task of identifying and learning such a representation is crucial to model reduction theory and methods. 

Traditionally, linear approaches such as the Proper Orthogonal Decomposition (POD) are the method of choice in problems with high-dimensional state spaces associated with physics-based data and modeling.\cite{lumley1967structures, sirovich1987turbulence, holmes1996turbulence} 
There are several model-reduction approaches based on POD, such as dynamic mode decomposition, balanced POD, the reduced basis method, POD-based discrete empirical interpolation, and data-driven operator inference.\cite{schmid2010dynamic, tu2013dynamic, kutz2016dynamic, doi:10.1146/annurev-fluid-030121-015835, willcox2002balanced, rowley2005model, quarteroni2015reduced, veroy2005certified, doi:10.1137/090766498, PEHERSTORFER2016196, QIAN2020132401} 
While all such methods have their own strengths and limitations, their formulation in terms of linear dimension reduction principles can lead to difficulties in problems for which the Kolmogorov $N$-width decreases slowly with increasing $N$. The Kolmogorov $N$-width is a measure for the worst-case error that might arise from the projection of a solution manifold onto a linear subspace of dimension $N$.\cite{Kolmogoroff} 

Nonlinear methods with machine learning techniques at their core have received have received a surge in attention over the past few years for adoption in reduced-order modeling applications. \cite{LEE2020108973, wan2018data, doi:10.1137/18M1177263, KIM2022110841, FRESCA2022114181} We here focus on methods that build on the concept of nonlinear \emph{embeddings}.\cite{JAIN201780, BARNETT2022111348, GEELEN2023115717, axaas2022fast, geelen2023learning, kalur2023data} 
More specifically, this paper develops a methodology in which we learn a reduced, nonlinear representation of the data, after which the system dynamics can be expressed and learned in terms of this new representation. 
We first reveal the low-dimensional manifold structure of a dynamical system by fitting a nonlinear basis expansion to the observed data. 
(The task of characterizing and representing this inherent structure is often referred to as {\em representation learning} in the machine learning community.\cite{6472238})
The construction of nonlinear state representations is achieved here by enriching linear approximations with low-order polynomial terms. 
We then summarize two different approaches for learning the manifold geometry from data.
The first is based on the POD, while the second builds on alternating minimization techniques.\cite{geelen2023learning} 
Once the nonlinear manifold can be represented in the form of a basis expansion, a projection of the PDE model onto the manifold reveals the algebraic structure of the reduced-order model for the problem of interest. 
Data-driven operator inference (OpInf) is then used to construct projection-based reduced-order models from data in a non-intrusive fashion.\cite{PEHERSTORFER2016196, QIAN2020132401} 

A major benefit of our proposed methodology is that the system states are expressed through a compact \emph{modal} representation. The use of modal decomposition and analysis techniques has paved the way to scalable nonlinear model reduction methods that favor interpretability and promote physical intuition. Modal representations have a long history in approximations of physics-based models, stretching back to early work in normal mode representations in structures and reduced basis representations. \cite{doi:10.2514/3.4741, KERSCHEN2009170, doi:10.2514/3.50778}. In many applications the modes have physical significance in that they represent dominant physical modal responses. \cite{dowell1996eigenmode, 10.1115/1.3101718, rowley_mezic_bagheri_schlatter_henningson_2009, doi:10.1146/annurev-fluid-011212-140652} Our proposed methodology follows this template in that the latent space coordinates are the multipliers (i.e., modal coordinates) of a nonlinear combination of basis vectors with physical interpretation. An immediate consequence of this construction is that the coordinates now also represent the \emph{evolution} of the dynamical system on a low-dimensional manifold. Because the proposed techniques are directly applicable to discrete-time observed data, the reduced model constructed in this fashion has the form of a system of nonlinear ordinary differential equation (ODEs), which can be solved effectively by modern ODE solvers.

The proposed approach shares conceptual parallels with the seminal work from Kevrekidis and coworkers on so-called approximate inertial manifolds.\cite{foias1988computation, JOLLY199038, johnson1997two, graham1996alternative} The theory of inertial manifolds states that infinite-dimensional system of partial different equations (PDEs) may be described accurately in their long term behavior by finite-dimensional systems. The existence of inertial manifolds has been established for various PDE systems in computational physics.\cite{foias1988computation, constantin1989integral, nicolaenko1989some, mallet1988inertial}
Approximate inertial manifold calculations are carried out as follows. One first performs a Galerkin approximation in the derivation of low-dimensional reduced systems, after which the higher-order modes are represented by means of the lower-order modes (a process sometimes referred to as slaving). 
Of particular relevance to this paper is the idea suggested in Dean et al.\cite{deane1991low}: "\emph{A particularly interesting direction is to combine the approximate inertial manifold method with the POD eigenmode hierarchy, and approximate the solution component on the higher POD modes as a function of its components on the lower, more energetic ones}".

An outline of the remainder of the paper follows.
Section~\ref{sec:learning_nonlinear_manifolds} discusses the construction of nonlinear state approximations through the lens of representation learning. 
We demonstrate carefully how conventional linear state approximations can be enriched with polynomial terms and how the unknown basis matrices, coefficient matrix, and  representation of the data in a low-dimensional coordinate system can be determined in a principled manner. 
Section~\ref{sec:roms} works with these nonlinear state approximations to derive the algebraic structure of the corresponding reduced-order models and shows how physics-based reduced-order models may be learned from data using the OpInf methodology. 
We then provide numerical evidence for the effectiveness of manifold-based OpInf approaches in Section~\ref{sec:numerical_experiments} using representative numerical experiments. 
Some conclusions and future research directions are presented in Section~\ref{sec:conclusions}.

\section{Learning nonlinear manifolds}
\label{sec:learning_nonlinear_manifolds}

This section outlines our method for learning nonlinear manifolds and presents an illustrative example.
Section~\ref{subsec:polynomial_manifolds} discusses the general representation learning problem for constructing nonlinear state approximations in problems with high-dimensional state spaces. 
The numerical procedures introduced by Geelen et al.\cite{geelen2023learning} for solving this learning problem are summarized in Section~\ref{subsec:basis_computation}.
In Section~\ref{subsec:illustrative_example}, a geometric interpretation of the method is illustrated by means of a simple three-dimensional example. 
A Jupyter notebook outlining the computational steps for this example is available at \url{https://github.com/geelenr/nl_manifolds}.

\subsection{A general representation learning problem}
\label{subsec:polynomial_manifolds}

Our focus is on data generated from complex PDE models that  represent the governing laws of nature. 
A training data set is comprised of a set of snapshots, each snapshot being a sample of the high-fidelity state representing a particular condition of the physical system. 
We denote each snapshot by $\mathbf{s}_j \in \mathbb{R}^n$ for $n \in \mathbb{N}$, and construct a snapshot matrix $\mathbf{S} \in \mathbb{R}^{n \times k}$ from $k$ such snapshots:
\begin{equation}
        \mathbf{S} := 
    \begin{pmatrix}
        | & | & & | \\
        \mathbf{s}_1 & \mathbf{s}_2  & \dots & \mathbf{s}_k  \\
        | & | & & |
    \end{pmatrix}.
    \label{eq:snapshot_data}
\end{equation}
We make use of a reference state, $\mathbf{s}_\text{ref} \in \mathbb{R}^n$, and denote by $\mathbf{S}_\text{ref}$ the $n \times k$ reference matrix each of whose columns is  $\mathbf{s}_\text{ref}$. 

To approximate the high-dimensional state $\mathbf{s}(t)\in \mathbb{R}^n$, we seek low-dimensional approximations $\boldsymbol{\Gamma}: \mathbb{R}^r \mapsto \mathbb{R}^n$ such that 
\begin{equation}
    \mathbf{s}(t) \approx \boldsymbol{\Gamma}(\hat{\mathbf{s}}(t)),
    \label{eq:nonlinear_mapping}
\end{equation}
where $t \in \mathbb{R}$ is some parameter on which the state depends (frequently, time). 
The vector $\hat{\mathbf{s}}(t) \in \mathbb{R}^r$ denotes the reduced state coordinate vector of dimension $r$. 
The transformation $\boldsymbol{\Gamma}$ constitutes a nonlinear mapping from the reduced-state coordinate system to the original, high-dimensional state space. 
We consider the following specific nonlinear modal basis structure for $\boldsymbol{\Gamma}$:
\begin{equation}
    \mathbf{s}(t) \approx \boldsymbol{\Gamma}(\hat{\mathbf{s}}(t)) := \mathbf{s}_\text{ref} + \underbrace{\mathbf{V} \hat{\mathbf{s}}(t)}_\text{linear} + \underbrace{\overline{\mathbf{V}} \boldsymbol{\Xi} \mathbf{g}(\hat{\mathbf{s}}(t))}_\text{nonlinear},
    \label{eq:nonlinear_approx}
\end{equation}
where $\mathbf{V}  = [\mathbf{v}_1 \,| \, \dotsc \,|\, \mathbf{v}_r] \in \mathbb{R}^{n \times r}$ and $\overline{\mathbf{V}} = [\overline{\mathbf{v}}_1 \,| \, \dotsc \,|\, \overline{\mathbf{v}}_{q}] \in \mathbb{R}^{n \times q}$ are a pair of basis matrices. The matrix $\boldsymbol{\Xi} \in \mathbb{R}^{q \times (p-1)r}$ is a coefficient matrix that controls the weighting of the basis functions contained in $\overline{\mathbf{V}}$. The vector $\mathbf{g}(\hat{\mathbf{s}}(t)) \in \mathbb{R}^{(p-1)r}$ has the form
\begin{equation}
    \mathbf{g}(\hat{\mathbf{s}}(t)) = \begin{pmatrix} \hat{\mathbf{s}}^2(t) \\
    \hat{\mathbf{s}}^3(t) \\
    \vdots                \\
    \hat{\mathbf{s}}^p(t) 
    \end{pmatrix},
    \label{eq:polynomial}
\end{equation}
where each $\hat{\mathbf{s}}^j(t) \in \mathbb{R}^r$ consists of the $j$th power of the components of $\hat{\mathbf{s}}(t)$, that is, $\hat{\mathbf{s}}^j(t) = [\hat{s}_1(t)^j,\hat{s}_2(t)^j,\dotsc,\hat{s}_r(t)^j]$. While one could instead employ Kronecker products of the reduced-state vector $\hat{\mathbf{s}}(t)$ (in contrast to the proposed element-wise formulation without the cross terms), the number of terms in $\mathbf{g}(\hat{\mathbf{s}}(t))$ would grow exponentially in such a construction. Keeping the dimension of $\mathbf{g}(\hat{\mathbf{s}}(t))$ relatively small is of particular importance in learning physics-based reduced-order models, to be discussed in Section \ref{sec:roms}.

The representation learning problem to construct a nonlinear state approximation of the form \eqref{eq:nonlinear_approx} is now posed as a constrained optimization problem\cite{geelen2023learning}
\begin{equation}
\begin{aligned}
    & \min_{\mathbf{V}, \overline{\mathbf{V}}, \boldsymbol{\Xi}, \hat{\mathbf{S}}} \left( F(\mathbf{V}, \overline{\mathbf{V}}, \boldsymbol{\Xi}, \hat{\mathbf{S}}) + \dfrac{\gamma }{2} \left\| \boldsymbol{\Xi} \right\|_F^2 \right) \\
    & \text{such that } \begin{pmatrix} \mathbf{V} & \overline{\mathbf{V}} \end{pmatrix} \in \mathcal{V}_{n,(r+q)},
    \label{eq:optim_problem_reg}
\end{aligned}
\end{equation}
where $ \hat{\mathbf{S}} := ( \hat{\mathbf{s}}_1 , \hat{\mathbf{s}}_2 , \dots, \hat{\mathbf{s}}_k ) \in \mathbb{R}^{r \times k} $ is the reduced-state representation of the given system states $\mathbf{s}_j$ for $j=1,\dots,k$; the objective function term $F$ is defined as 
\begin{equation}
\begin{aligned}
    F(\mathbf{V}, \overline{\mathbf{V}}, \boldsymbol{\Xi}, \hat{\mathbf{S}}) &= \dfrac{1}{2} \sum_{j=1}^k  \left\| \mathbf{s}_j - \boldsymbol{\Gamma}(\hat{\mathbf{s}}_j) \right\|_2^2 \\
    &= \dfrac{1}{2} \sum_{j=1}^k  \left\| \mathbf{s}_j - \mathbf{s}_\text{ref} - \begin{pmatrix} \mathbf{V} & \overline{\mathbf{V}} \end{pmatrix} \begin{pmatrix} \hat{\mathbf{s}}_j \\ \boldsymbol{\Xi} \mathbf{g}( \hat{\mathbf{s}}_j ) \end{pmatrix} \right\|_2^2;
    \label{eq:objective}
\end{aligned}
\end{equation}
Frobenius norm regularization involving $\boldsymbol{\Xi}$ is used to to avoid overfitting to the training data; $\gamma \ge 0$ is a regularization parameter; and $\mathcal{V}_{n,(r+q)}$ is the Stiefel manifold, defined as the set of matrices in $\mathbb{R}^{n \times (r+q)}$ with orthonormal columns, that is,
\begin{equation}
    \mathcal{V}_{n,(r+q)} = \{ \begin{pmatrix} \mathbf{V} & \overline{\mathbf{V}} \end{pmatrix} \in \mathbb{R}^{n \times (r+q)}: \begin{pmatrix} \mathbf{V} & \overline{\mathbf{V}} \end{pmatrix}^\top \begin{pmatrix} \mathbf{V} & \overline{\mathbf{V}} \end{pmatrix} = \mathbf{I}_{r+q} \},
    \label{eq:orthogonality}
\end{equation} 
where $\mathbf{I}_{r+q}$ is the $\mathbb{R}^{(r+q) \times (r+q)}$ identity matrix.

It should be noted that an expressive approximation of the form \eqref{eq:nonlinear_approx} has the tendency to overfit noise or anomalous behavior in the data. 
In general, the more parameters in the model, the higher the likelihood that it will overfit. Specifically, approximations with only low-order polynomial degree $p$ (as we advocate here) are  less prone to overfitting in comparison to higher order polynomial methods or alternative black-box regression methods. 
The addition of Frobenius regularization in \eqref{eq:optim_problem_reg} also helps to mitigate overfitting.

\subsection{Computing the basis expansion}
\label{subsec:basis_computation}

\begin{figure}[!tbp]
\begin{algorithm}[H]
\caption{POD-based representation learning}
\label{alg:pod-based}
\begin{algorithmic}
\REQUIRE Snapshot matrix $\mathbf{S} \in \mathbb{R}^{n \times k}$, reference state $\mathbf{s}_\text{ref} \in \mathbb{R}^n$, regularization $\gamma \in \mathbb{R}^+$, polynomial order $p \in \mathbb{N}_{\geq 2}$.
\ENSURE Basis matrices $\mathbf{V} \in \mathbb{R}^{n \times r}$ and $ \overline{\mathbf{V}} \in \mathbb{R}^{n \times q}$, coefficient matrix $\boldsymbol{\Xi} \in \mathbb{R}^{q \times (p-1)r}$, projected snapshot data $\hat{\mathbf{S}} \in \mathbb{R}^{r \times k}$.
\STATE \textbf{Step 1}: Fixing an orthogonal set of basis vectors
\begin{equation}
    \mathbf{S}-\mathbf{S}_\text{ref}=\boldsymbol{\Phi} \boldsymbol{\Sigma} \boldsymbol{\Psi}^\top.
\end{equation}
\STATEx $(\mathbf{V},\overline{\mathbf{V}}) \leftarrow$ the $r+q$ leading left singular vectors of $\mathbf{S}-\mathbf{S}_\text{ref}$
\STATE \textbf{Step 2}: Compute snapshot representation in POD coordinates
\begin{equation}
\hat{\mathbf{S}} = \mathbf{V}^\top (\mathbf{S} - \mathbf{S}_\text{ref}).
\end{equation}
\STATE \textbf{Step 3}: Fixing $\mathbf{V}$, $\overline{\mathbf{V}}$ and $\hat{\mathbf{S}}$ at their values defined above, compute the coefficient matrix $\boldsymbol{\Xi}$ by solving a linear least-squares problem:
\begin{equation}
    \min_{\boldsymbol{\Xi}} \left( F(\mathbf{V}, \overline{\mathbf{V}}, \boldsymbol{\Xi}, \hat{\mathbf{S}}) + \dfrac{\gamma}{2} \left\| \boldsymbol{\Xi} \right\|_F^2 \right).
    \label{eq:solving_for_psi}
\end{equation}
\end{algorithmic}
\end{algorithm}
\end{figure}

\begin{figure}[!tbp]
\begin{algorithm}[H]
\caption{Alternating minimization based representation learning}
\label{alg:am-based}
\begin{algorithmic}
\REQUIRE Snapshot matrix $\mathbf{S} \in \mathbb{R}^{n \times k}$, reference state $\mathbf{s}_\text{ref} \in \mathbb{R}^n$, regularization $\gamma \in \mathbb{R}^+$, polynomial order $p \in \mathbb{N}_{\geq 2}$, stopping criterion.
\ENSURE Basis matrices $\mathbf{V} \in \mathbb{R}^{n \times r}$ and $\overline{\mathbf{V}} \in \mathbb{R}^{n \times q}$, coefficient matrix $\boldsymbol{\Xi} \in \mathbb{R}^{q \times (p-1)r}$, projected snapshot data $\hat{\mathbf{S}} \in \mathbb{R}^{r \times k}$.
\WHILE {stopping criterion not satisfied}
\STATE \textbf{Step 1}: Orthogonal Procrustes: Compute the basis vectors $\boldsymbol{\Omega} := \begin{pmatrix} \mathbf{V} & \overline{\mathbf{V}} \end{pmatrix}$ by solving
\begin{equation}
\begin{aligned}
    &\min_{\boldsymbol{\Omega}} \, \dfrac{1}{2} \left\| \mathbf{S} - \mathbf{S}_\text{ref} - \boldsymbol{\Omega} \begin{pmatrix} \hat{\mathbf{S}} \\ \boldsymbol{\Xi} \mathbf{g}(\hat{\mathbf{S}}) \end{pmatrix} \right\|_F^2 \\ &\text{such that } \boldsymbol{\Omega}^\top \boldsymbol{\Omega} = \mathbf{I}_{r+q}.
    \label{eq:procrustes}
\end{aligned}
\end{equation}
\STATE \textbf{Step 2}: Compute the coefficient matrix $\boldsymbol{\Xi}$ from a linear least-squares problem:
\begin{equation}
    \min_{\boldsymbol{\Xi}} \, \left( F(\mathbf{V}, \overline{\mathbf{V}}, \boldsymbol{\Xi}, \hat{\mathbf{S}}) + \dfrac{\gamma}{2} \left\| \boldsymbol{\Xi} \right\|_F^2 \right)
    \label{eq:solving_for_psi_2}
\end{equation}
\STATE \textbf{Step 3}: Project data  onto the nonlinear manifold by solving
\begin{equation}
    \min_{\hat{\mathbf{s}}_j} \, \dfrac{1}{2}  \left\| \mathbf{s}_j - \mathbf{s}_\text{ref} - \boldsymbol{\Omega} \begin{pmatrix} \hat{\mathbf{s}}_j \\ \boldsymbol{\Xi} \mathbf{g}( \hat{\mathbf{s}}_j ) \end{pmatrix} \right\|_2^2, \;\; j=1,2,\dotsc,k. 
    \label{eq:opt_low_dim_coord_system}
\end{equation}
\ENDWHILE
\end{algorithmic}
\end{algorithm}
\end{figure}

This section summarizes the two methodologies from Geelen et al.\cite{geelen2023learning} for finding a numerical approximation to the solution of the  representation learning problem \eqref{eq:optim_problem_reg}--\eqref{eq:objective}: the POD-based and alternating-minimization-based methods. 
To build approximations of the form \eqref{eq:nonlinear_approx}, these approaches make informed choices on the basis matrices $\mathbf{V}$ and $\overline{\mathbf{V}}$, the coefficient matrix $\boldsymbol{\Xi}$, and the reduced-state representation of the data $\hat{\mathbf{S}}$.

In the POD-based representation learning method (see  Algorithm~\ref{alg:pod-based}), the columns of $\mathbf{V}$ are chosen to be the POD basis vectors, that is, the left singular vectors of $\mathbf{S}-\mathbf{S}_\text{ref}$ corresponding to the $r$ largest singular values. 
The columns of the basis matrix $\overline{\mathbf{V}}$ are chosen to be the left singular vectors corresponding to the next $q$ largest singular values. 
This choice of matrices \{$\mathbf{V}, \overline{\mathbf{V}}$\} satisfies constraint \eqref{eq:orthogonality} by virtue of the orthogonality property of singular vectors. 
Representation of snapshots in the POD coordinates can be calculated by means of an orthogonal projection. 
The coefficient matrix $\boldsymbol{\Xi}$ is obtained from a linear least-squares problem. 
% The method, which we call POD-based representation learning, is summarized in Algorithm~\ref{alg:pod-based}), 

\begin{figure*}[!tbp]
\centering \footnotesize
\begin{subfigure}[t]{.32\linewidth}
\includegraphics[width=\linewidth]{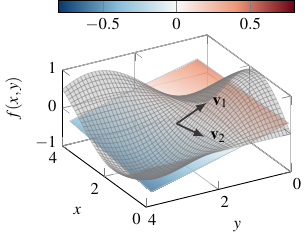}
  \caption{Traditional POD (linear subspace) ($22.6\%$ error)}
\end{subfigure}
\begin{subfigure}[t]{.32\linewidth}
\includegraphics[width=\linewidth]{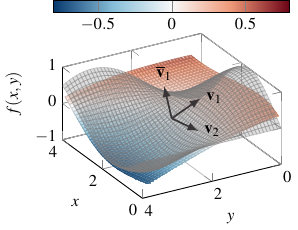}
  \caption{POD-based representation learning ($21.1\%$ error)}
  \label{fig:3d_pod_manifold}
\end{subfigure} 
\begin{subfigure}[t]{.32\linewidth}
\includegraphics[width=\linewidth]{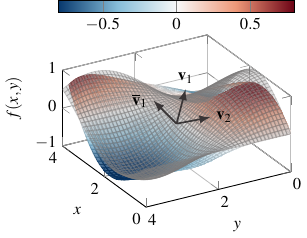}
  \caption{AM-based representation learning ($10.4\%$ error)}
  \label{fig:3d_am_manifold}
\end{subfigure}
\caption{Comparison of (a) the linear-subspace POD with (b) the POD-based representation learning and (c) alternating minimization based representation learning approaches for reconstructing, using two modal coefficients, the three-dimensional manifold. The relative state error $\| \mathbf{S}-\boldsymbol{\Gamma}(\hat{\mathbf{S}}) \|_F/ \| \mathbf{S} - \mathbf{S}_\text{ref} \|_F$ is given between parentheses. The gray surface denotes the original three-dimensional manifold, whereas the colored surfaces illustrate the different reconstructions. The black arrows represent the basis vectors.}
\label{fig:3d_problem}
\end{figure*}

The alternating minimization (AM) representation learning approach (Algorithm~\ref{alg:am-based}) proceeds as follows. 
Initial guesses are required for the projected snapshot data matrix $\hat{\mathbf{S}}$ and the coefficient matrix $\boldsymbol{\Xi}$. 
Then the objective in \eqref{eq:optim_problem_reg} is successively minimized for the three blocks of variables in turn --- first $\boldsymbol{\Omega} = \begin{pmatrix} \mathbf{V} & \overline{\mathbf{V}} \end{pmatrix}$, then $\boldsymbol{\Xi}$, then $\hat{\mathbf{S}}$ --- with the pattern repeating until  a convergence criterion is satisfied. 
Generally speaking, three-block alternating minimization schemes have no theoretical guarantees of convergence. \cite{chen2016direct} However in this context, and many others, they converge on practical instances. \cite{geelen2023learning}

The minimization  with respect to $\boldsymbol{\Omega}$ is a standard problem known as the orthogonal Procrustes problem, and it can be solved via singular value decomposition. \cite{schonemann1966generalized}
The minimization with respect to $\boldsymbol{\Xi}$ is a linear least squares problem, as in Algorithm~\ref{alg:pod-based}.
The minimization with respect to $\hat{\mathbf{S}}$ decomposes into $k$ separate problems, each of which has a single reduced state $\hat{\mathbf{s}}_j$ as its variable. 

\subsection{Orthogonal subspace transformations: an illustrative example}
\label{subsec:illustrative_example}

We demonstrate the nonlinear representation \eqref{eq:nonlinear_approx}  by means of a small numerical example.
Consider a manifold in a three-dimensional Euclidian space parametrized by the vector $\mathbf{s}(x,y) = \left( x, y, \sin(x)\cos(y) \right)^\top \in \mathbb{R}^3$ for $x, y \in [0,4]$. 
A dataset is built by sampling $\mathbf{s}(x,y)$ uniformly over the domain with grid spacings $\Delta x=\Delta y = 0.1$. 
The resulting data matrix $\mathbf{S}$ has dimension $3\times 1,681$. 
After defining $\mathbf{s}_\text{ref}$ to be the column-averaged mean of the data matrix, we build approximations to the nonlinear manifold of dimension $r=2$. 
By choosing a polynomial embedding of degree $p=3$ in \eqref{eq:polynomial}, the nonlinear state approximation \eqref{eq:nonlinear_approx} of the $j$th data sample becomes
\begin{equation}
\begin{aligned}
    \mathbf{s}_j &\approx \mathbf{s}_\text{ref} + \underbrace{\mathbf{v}_1 \hat{s}_{1,j} + \mathbf{v}_2 \hat{s}_{2,j}}_\text{linear} \\
    &\quad + \underbrace{\overline{\mathbf{v}}_1 \left( \Xi_1 \hat{s}_{1,j}^2 + \Xi_2 \hat{s}_{2,j}^2 + \Xi_3 \hat{s}_{2,j}^3 + \Xi_4 \hat{s}_{2,j}^3 \right)}_\text{nonlinear},
    \label{eq:3d_nonlinear}
\end{aligned}
\end{equation}
where the basis vectors $\mathbf{v}_1, \mathbf{v}_2, \overline{\mathbf{v}}_1$ form an orthogonal set. The modal coefficients are given by $\hat{s}_{1,j}$ and $\hat{s}_{2,j}$. Instead of introducing a \emph{separate} coefficient for the third basis vector $\overline{\mathbf{v}}_1$, we express its coefficient in terms of the coefficients of the first two basis vectors $\hat{s}_{1,j}$ and $\hat{s}_{2,j}$. 
The coefficients $\boldsymbol{\Xi} = (\Xi_1,\Xi_2,\Xi_3,\Xi_4) \in \mathbb{R}^{4}$ control the weighting of $\overline{\mathbf{v}}_1$.

We now follow the steps from the POD-based representation learning formulation of Algorithm~\ref{alg:pod-based}. 
The vectors $\mathbf{v}_1, \mathbf{v}_2, \overline{\mathbf{v}}_1$ are fixed to be the three left singular vectors of the shifted data matrix. 
Accordingly, $\hat{s}_{1,j}$ and $\hat{s}_{2,j}$ represent the coefficients of expansion in the basis $\mathbf{V}$. 
The coefficients $\boldsymbol{\Xi}$ are inferred from the data via linear regression; see \eqref{eq:solving_for_psi}.

In the alternating minimization based representation learning approach of Algorithm~\ref{alg:am-based}, the basis vectors $\mathbf{v}_1, \mathbf{v}_2, \overline{\mathbf{v}}_1$ and the data representation in the reduced-space coordinate system (through coordinates $\hat{s}_1, \hat{s}_2$) are computed by way of an orthogonal Procrustes problem and a set of  $k=1,681$ unconstrained nonlinear optimization problems, respectively. 
The computation of the coefficient matrix $\boldsymbol{\Xi}$ is the linear least-squares problem \eqref{eq:solving_for_psi_2}. 

Fig.~\ref{fig:3d_problem} compares the reconstructions of the standard (linear) POD approximation with those that use \eqref{eq:3d_nonlinear}. 
Linear-subspace POD invokes a large projection error and, in this example, is ill-suited for data reconstruction tasks. 
The POD-based representation learning approach warps the POD subspace to produce a nonlinear manifold that is slightly closer to the exact solution, see Fig.~\ref{fig:3d_pod_manifold}. 
By applying a series of rotations and/or reflections to the POD basis, the alternating minimization based method finds preferred directions along which to apply curvature to further reduce the representation error, as shown in Fig.~\ref{fig:3d_am_manifold}.

\section{Learning reduced-order models on nonlinear manifolds}
\label{sec:roms}

In this section, we show how nonlinear manifold representations of the form \eqref{eq:nonlinear_approx} can be employed to learn physics-based reduced-order model from data. 
Assuming that the full-space data $\mathbf{s}(t)$ is generated by solving nonlinear governing equations with particular (widely relevant) structure, we substitute from \eqref{eq:nonlinear_approx} to obtain the corresponding system in the reduced space. We then describe a process for learning the operators that define the reduced-order model.

Section~\ref{subsec:projection-based_roms} derives the algebraic structure of the underlying reduced-order models through a manifold projection method. 
In Section~\ref{subsec:opinf}, we propose a manifold-based inference method for constructing reduced-order models from snapshot data.

\subsection{Projection-based model reduction}
\label{subsec:projection-based_roms}

Consider the following initial-value nonlinear ODE problem:
\begin{equation}
    \dfrac{\text{d}}{\text{d}t} \mathbf{s}(t) =  \mathbf{f}(\mathbf{s}(t)); \quad \mathbf{s}(0) = \mathbf{s}_0,
    \label{eq:nonlinear_pde}
\end{equation}
where as before $\mathbf{s}(t) \in \mathbb{R}^n$ is the system state at time $t$ and $\mathbf{f}: \mathbb{R}^n \mapsto \mathbb{R}^n$ maps the state to its time derivative.
% , with initial condition is given by $\mathbf{s}(0) = \mathbf{s}_0$. 
In many systems that arise throughout computational engineering and sciences, the operator $\mathbf{f}$ has a certain linear-quadratic form, \cite{5991229, doi:10.2514/1.J057791, QIAN2020132401} allowing us to work with the following special case of \eqref{eq:nonlinear_pde}:
\begin{equation}
\dfrac{\text{d}\mathbf{s} }{\text{d}t} =  \mathbf{A}\mathbf{s} + \mathbf{H}(\mathbf{s} \otimes \mathbf{s}); \quad \mathbf{s}(0) = \mathbf{s}_0,
\label{eq:fom}
\end{equation}
where we omitted the dependence of the system states on $t$ to simplify the notation, and $\otimes$ denotes the Kronecker product. 
(We will refer to  \eqref{eq:fom} henceforth as the full-order model (FOM).) 
The operators $\mathbf{A} \in \mathbb{R}^{n \times n}$ and $\mathbf{H} \in \mathbb{R}^{n \times n^2}$ denote the FOM operators corresponding to linear and quadratic terms, respectively, in the governing semi-discrete equations. 

The use of nonlinear state approximations of the form \eqref{eq:nonlinear_approx} informs the algebraic structure of the reduced-order analog  of \eqref{eq:fom}, requiring that we account for cubic and higher-order interactions between the modal coefficients. 
Specifically, by introducing  \eqref{eq:nonlinear_approx} into \eqref{eq:fom}, and projecting the residual of the resulting system onto the span of $\mathbf{V}$, we obtain
\begin{equation}
    \begin{aligned}
\dfrac{\text{d}\hat{\mathbf{s}}}{\text{d}t} &= \mathbf{V}^\top \mathbf{A} \mathbf{s}_\text{ref} + \mathbf{V}^\top \mathbf{A} \mathbf{V} \hat{\mathbf{s}} + \mathbf{V}^\top \mathbf{A} \overline{\mathbf{V}} \boldsymbol{\Xi} \mathbf{g}(\hat{\mathbf{s}}) + \mathbf{V}^\top \mathbf{H}(\mathbf{s}_\text{ref} \otimes \mathbf{s}_\text{ref}) \\ & \quad + \mathbf{V}^\top \mathbf{H}(\mathbf{s}_\text{ref} \otimes \mathbf{V}\hat{\mathbf{s}}) + \mathbf{V}^\top \mathbf{H}(\mathbf{s}_\text{ref} \otimes \overline{\mathbf{V}}\boldsymbol{\Xi}\mathbf{g}(\hat{\mathbf{s}})) \\ & \quad + \mathbf{V}^\top \mathbf{H}( \mathbf{V}\hat{\mathbf{s}} \otimes \mathbf{s}_\text{ref}) + \mathbf{V}^\top \mathbf{H}(  \mathbf{V}\hat{\mathbf{s}} \otimes \mathbf{V}\hat{\mathbf{s}}) \\ & \quad + \mathbf{V}^\top \mathbf{H}(  \mathbf{V}\hat{\mathbf{s}} \otimes \overline{\mathbf{V}}\boldsymbol{\Xi}\mathbf{g}(\hat{\mathbf{s}}))  + \mathbf{V}^\top \mathbf{H}( \overline{\mathbf{V}}\boldsymbol{\Xi}\mathbf{g}(\hat{\mathbf{s}}) \otimes \mathbf{s}_\text{ref}) \\
& \quad + \mathbf{V}^\top \mathbf{H}(  \overline{\mathbf{V}}\boldsymbol{\Xi}\mathbf{g}(\hat{\mathbf{s}}) \otimes \mathbf{V}\hat{\mathbf{s}}) + \mathbf{V}^\top \mathbf{H}(  \overline{\mathbf{V}}\boldsymbol{\Xi}\mathbf{g}(\hat{\mathbf{s}}) \otimes \overline{\mathbf{V}}\boldsymbol{\Xi}\mathbf{g}(\hat{\mathbf{s}})),
\label{eq:projection-based-rom}
    \end{aligned}
\end{equation}
with $\hat{\mathbf{s}}(0) = \hat{\mathbf{s}}_0$ the representation of the initial condition $\mathbf{s}_0$ in the low-dimensional coordinate system. An alternative approach would be to choose a nonlinear, state-dependent projection of the state equations, which comes at the cost of increased algebraic complexity in the reduced-order models\cite{benner2023quadratic}.

By adapting the argument of Section~2.2 from a previous work\cite{GEELEN2023115717} to this case, we can show that when $\mathbf{V}$ and $\hat{\mathbf{S}}$ are obtained from a SVD of $\mathbf{S}-\mathbf{S}_\text{ref}$ (as in standard linear POD), $\boldsymbol{\Xi}$ is fixed at some reasonable value, and we minimize the objective function in \eqref{eq:optim_problem_reg} over $\overline{\mathbf{V}}$ alone, the orthogonality property $\mathbf{V}^\top \overline{\mathbf{V}}=\mathbf{0}$ will be satisfied automatically by the minimizing value of $\overline{\mathbf{V}}$. Although this property cannot be guaranteed to hold when we optimize {\em simultaneously} with respect to $\mathbf{V}$, $\overline{\mathbf{V}}$, $\hat{\mathbf{S}}$, and $\boldsymbol{\Xi}$, as in \eqref{eq:optim_problem_reg}, much remains to be gained in terms of simplification of our inferred reduced-order model. Specifically, by enforcing the orthogonality property in \eqref{eq:optim_problem_reg} the left-hand side of \eqref{eq:projection-based-rom} simplifies to a time derivative applied to the reduced state vector $\hat{\mathbf{s}}(t)$. A lack of orthogonality, on the other hand, would result in an expression that involves the time derivative of the polynomial function $\mathbf{g}(\hat{\mathbf{s}}(t))$. The added nonlinearity would cause the evaluation of the data-driven models to be more mathematically cumbersome.

The right-hand side of \eqref{eq:projection-based-rom} contains polynomial nonlinear terms up to order $2p$. In practice this structure poses a significant challenge from the implementation viewpoint: explicitly computing the different projected operators becomes cumbersome and requires explicit access to the full-order operators $\mathbf{A}$ and $ \mathbf{H}$. 
Rather than operating on \eqref{eq:projection-based-rom} directly, we expose its polynomial structure by using the mixed-product property of Kronecker products and grouping the constant, linear, quadratic, and higher-order terms as follows: 
\begin{equation}
\dfrac{\text{d}\hat{\mathbf{s}}}{\text{d}t} = \hat{\mathbf{c}} + \hat{\mathbf{A}} \hat{\mathbf{s}} + \hat{\mathbf{H}} (\hat{\mathbf{s}} \otimes \hat{\mathbf{s}})  + \hat{\mathbf{P}}\hat{\mathbf{g}}(\hat{\mathbf{s}})
\label{eq:rom_nonlinear}
\end{equation}
where $\hat{\mathbf{c}} \in \mathbb{R}^r, \hat{\mathbf{A}} \in \mathbb{R}^{r \times r}, \hat{\mathbf{H}} \in \mathbb{R}^{r \times r^2}, \hat{\mathbf{P}} \in \mathbb{R}^{r \times d(r,p)}$ are the reduced matrix operators. 
The operator $\hat{\mathbf{P}}$ accounts for the higher-order interactions between the modal coefficients in the reduced-order model. 
These interactions are captured by the vector $\hat{\mathbf{g}}(\hat{\mathbf{s}})$ (which is a subvector of $\mathbf{g}(\hat{\mathbf{s}})$ defined in \eqref{eq:polynomial})  and consist of monomials from degree three to degree $2p$. 
The total number of unique coefficients in $\hat{\mathbf{g}}(\hat{\mathbf{s}})$ scales as $d(r,p) \sim \mathcal{O}(p^2r^2)$. 
Example~\ref{ex:example} provides an illustration of the structure of $\hat{\mathbf{g}}(\hat{\mathbf{s}})$ and the scaling of $d(r,p)$ in reduced-order models of the form \eqref{eq:rom_nonlinear}. 
A technique for approximating the reduced-order operators $\hat{\mathbf{c}}, \hat{\mathbf{A}}, \hat{\mathbf{H}}, \hat{\mathbf{P}}$ is presented in Section~\ref{subsec:opinf}.

\begin{exmp}\label{ex:example}
For illustration purposes, consider a state approximation of dimensionality $r=2$ and polynomial degree $p=3$ given by
\begin{equation}
    \mathbf{s} \approx \mathbf{s}_\text{ref} 
    + \mathbf{V} (\hat{s}_1, \hat{s}_2 )^\top 
    + \overline{\mathbf{V}} \boldsymbol{\Xi} 
    (\hat{s}_1^2, \hat{s}_2^2, \hat{s}_1^3, \hat{s}_2^3 )^\top,
    \label{eq:toy_problem}
\end{equation}
where $\mathbf{V} \in \mathbb{R}^{n \times 2}$ and $\overline{\mathbf{V}} \in \mathbb{R}^{n \times q}$ are the basis matrices and $\boldsymbol{\Xi} \in \mathbb{R}^{q \times 4}$ is the coefficient matrix, calculated as described in Section~\ref{sec:learning_nonlinear_manifolds}. 
The modal coefficients $\hat{s}_1, \hat{s}_2$ are the only unknowns. 
By substituting from \eqref{eq:toy_problem} into the linear-quadratic full-order model \eqref{eq:fom} and following the derivation from Section~\ref{subsec:projection-based_roms}, we obtain a reduced model of the form \eqref{eq:rom_nonlinear}. 
This model then accounts for the following nonlinear interactions among the modal coefficients: 
\begin{equation}
\begin{aligned}
    \hat{\mathbf{s}} &= \left( \hat{s}_1, \hat{s}_2\right)^\top \in \mathbb{R}^{2}, \\
    \hat{\mathbf{s}} \otimes \hat{\mathbf{s}} &= \left( \hat{s}_1^2, \hat{s}_1 \hat{s}_2, \hat{s}_2^2 \right)^\top \in \mathbb{R}^{3}, \\
    \hat{\mathbf{g}}(\hat{\mathbf{s}}) &= \left( \hat{s}_1^3, \hat{s}_1\hat{s}_2^2, \hat{s}_1^2\hat{s}_2, \hat{s}_2^3, \hat{s}_1^4, \hat{s}_1\hat{s}_2^3,  \hat{s}_1^3\hat{s}_2, \hat{s}_1^2\hat{s}_2^2, \hat{s}_2^4, \dots \right. \\ & \qquad \qquad \hat{s}_1^5, \hat{s}_1^2\hat{s}_1^3, \hat{s}_1^3\hat{s}_2^2, \hat{s}_2^5, \hat{s}_1^6, \hat{s}_1^3\hat{s}_2^3, \hat{s}_2^6  \left. \right)^\top \in \mathbb{R}^{d(r,p)=16},
\end{aligned}
\end{equation}
where $\hat{\mathbf{s}}, \hat{\mathbf{s}} \otimes \hat{\mathbf{s}}$, and $ \hat{\mathbf{g}}(\hat{\mathbf{s}})$ contain monomials of the modal coefficients of first, second, and higher-order degree, respectively. Using elementwise powers of the reduced-state vector $\hat{\mathbf{s}}$ (see \eqref{eq:polynomial}) ensures that the number of entries in $\hat{\mathbf{g}}(\hat{\mathbf{s}})$ remains tractable. Table \ref{tab:num_of_coefficients} lists $d(r,p)$ the number of terms contained in $\hat{\mathbf{s}}$ as a function of the reduced basis dimension $r$ and the degree of the polynomial embeddings $p$. 
\end{exmp}

\begin{table}[!tbp]
    \centering
    \caption{Dimension $d(r,p)$ as a function of the reduced basis dimension, $r$ and the degree of the polynomial embeddings, $p$.}
    \label{tab:num_of_coefficients}
    \begin{tabular}{|c|c|c|c|}
    \hline
         & $p=2$ & $p=3$ & $p=4$ \\
    \hline
        $r=2$  & 7   & 16  & 27  \\
        $r=4$  & 26  & 64  & 114 \\
        $r=6$  & 57  & 144 & 261 \\
        $r=8$  & 100 & 256 & 468 \\
        $r=10$ & 155 & 400 & 735 \\
    \hline
    \end{tabular}
\end{table}

\subsection{Learning physics-based reduced-order models from data}
\label{subsec:opinf}

We employ the data-driven operator inference (OpInf) method for learning the low-dimensional dynamical system \eqref{eq:rom_nonlinear} from time-domain simulation data.\cite{PEHERSTORFER2016196} 
While traditional linear-subspace POD lies at heart of the formulation of Peherstorfer and Willcox,\cite{PEHERSTORFER2016196} OpInf can be extended to the nonlinear manifold setting, as demonstrated in Geelen et al.\ for linear systems.\cite{GEELEN2023115717} Here, we consider reduction of nonlinear systems. The OpInf methodology finds the reduced matrix operators $\hat{\mathbf{c}}, \hat{\mathbf{A}}, \hat{\mathbf{H}}, \hat{\mathbf{P}}$ that define the reduced model that best matches the projected snapshot data, in the following sense of regularized least squares:
\begin{equation}
\begin{aligned}
\big( \hat{\mathbf{c}}, \hat{\mathbf{A}}, & \hat{\mathbf{H}}, \hat{\mathbf{P }} \big) = \argmin_{\tilde{\mathbf{c}}, \tilde{\mathbf{A}} , \tilde{\mathbf{H}}, \tilde{\mathbf{P }}} 
\Big( J(\tilde{\mathbf{c}}, \tilde{\mathbf{A}} , \tilde{\mathbf{H}}, \tilde{\mathbf{P }})  \\
&  + \dfrac{\lambda_1 }{2} \left( \| \tilde{\mathbf{c}} \|_2^2 + \| \tilde{\mathbf{A}} \|_F^2 \right)   + \dfrac{\lambda_2}{2} \| \tilde{\mathbf{H}} \|_F^2  + \dfrac{\lambda_3}{2} \| \tilde{\mathbf{P}} \|_F^2 \Big),
 \label{eq:opinf_nonlinear}
 \end{aligned}
\end{equation}
where the function $J(\tilde{\mathbf{c}}, \tilde{\mathbf{A}} , \tilde{\mathbf{H}}, \tilde{\mathbf{P }}) $ is defined to be
\begin{equation}
\sum_{j=1}^k
 \left\|  \tilde{\mathbf{c}} + \tilde{\mathbf{A}} \hat{\mathbf{s}}_j + \tilde{\mathbf{H}}(\hat{\mathbf{s}}_j \otimes \hat{\mathbf{s}}_j) + \tilde{\mathbf{P}}\hat{\mathbf{g}}(\hat{\mathbf{s}}_j) - \dfrac{\text{d}\hat{\mathbf{s}}_j}{\text{d}t} \right\|_2^2,
 \label{eq:opinf_objective}
\end{equation}
while the nonnegative scalars $\lambda_i$ with $i=1,2,3$ are Tikhonov regularization parameters that promote stability in the inferred reduced-order models and inhibit the overfitting of the system operators to potentially noisy data.\cite{doi:10.1080/03036758.2020.1863237} 
Generally speaking, good choices of the regularization parameters $\lambda_i$ tend to be different as they are coefficients of terms with different scales. 
The time derivatives in the objective function are typically estimated numerically using finite difference approximations. 
The optimization problem \eqref{eq:opinf_nonlinear} decouples into $r$ independent linear least-squares problems. \cite{PEHERSTORFER2016196} 

Algorithm~\ref{alg:standard_opinf} presents the steps of the linear-subspace OpInf approach for quadratic systems,\cite{PEHERSTORFER2016196, doi:10.1080/03036758.2020.1863237} while the workflow of the proposed nonlinear manifold-based OpInf methodology is summarized in Algorithm~\ref{alg:proposed_opinf}. 
We use the acronyms MPOD-OpInf and MAM-OpInf to distinguish between the nonlinear manifold OpInf approaches  based on Algorithm~\ref{alg:pod-based} and Algorithm~\ref{alg:am-based}, respectively.
These methodologies differ in the manner in which the projected snapshot data are computed and thus also in their subsequent reconstructions in the original state space. This is due to differences in the adopted low-dimensional basis, coefficient matrix, and reduced-state data representation; see Section~\ref{sec:learning_nonlinear_manifolds}.

\begin{figure}[!tbp]
\begin{algorithm}[H]
\caption{The standard (POD-based) OpInf methodology for quadratic
systems\cite{PEHERSTORFER2016196, doi:10.1080/03036758.2020.1863237}}
\label{alg:standard_opinf}
\begin{algorithmic}[1]
\REQUIRE Snapshot matrix $\mathbf{S} \in \mathbb{R}^{n \times k}$, reference state $\mathbf{s}_\text{ref} \in \mathbb{R}^n$
\ENSURE Reduced operators $\hat{\mathbf{c}}$, $\hat{\mathbf{A}}$, $\hat{\mathbf{H}}$
\STATE Compute SVD of $\mathbf{S}-\mathbf{S}_\text{ref}$
\STATE $\mathbf{V}$ $\leftarrow$ The $r$ leading left singular vectors of the $\mathbf{S}-\mathbf{S}_\text{ref}$
\STATE $\hat{\mathbf{S}}$ $\leftarrow$ $\mathbf{V}^\top (\mathbf{S}-\mathbf{S}_\text{ref})$ \hfill \Comment{Project snapshot data onto POD subspace}
\STATE Approximate $\dfrac{\text{d}}{\text{d}t}{\hat{\mathbf{S}}}$ from $\hat{\mathbf{S}}$ \hfill \Comment{Time derivative approximation}
\STATE Choose $\lambda_1 , \lambda_2$ \hfill \Comment{Hyperparameter optimization}
\STATE  Find reduced operators $\hat{\mathbf{c}}$, $\hat{\mathbf{A}}$, $\hat{\mathbf{H}}$ through OpInf regression
\end{algorithmic}
\end{algorithm}
\vspace{-1em}
\end{figure}

\begin{figure}[!tbp]
\begin{algorithm}[H]
\caption{Proposed nonlinear manifold based OpInf method for quadratic
systems}
\label{alg:proposed_opinf}
\begin{algorithmic}[1]
\REQUIRE Snapshot matrix $\mathbf{S} \in \mathbb{R}^{n \times k}$, reference state $\mathbf{s}_\text{ref} \in \mathbb{R}^n$, regularization parameter $\gamma \in \mathbb{R}^+$, polynomial order $p \in \mathbb{N}_{\geq 2}$, stopping criterion
\ENSURE Reduced operators $\hat{\mathbf{c}}$, $\hat{\mathbf{A}}$, $\hat{\mathbf{H}}$, $\hat{\mathbf{P}}$
\STATE Compute SVD of $\mathbf{S}-\mathbf{S}_\text{ref}$
\STATE $r,q$ $\leftarrow$ Choose number of left singular vectors to be used in basis matrices $\mathbf{V}$ and $\overline{\mathbf{V}}$
\STATE $\boldsymbol{\Gamma}(\mathbf{V}, \overline{\mathbf{V}}, \boldsymbol{\Xi}, \hat{\mathbf{S}})$ $\leftarrow$ Approximate the solution to general representation learning problem \eqref{eq:optim_problem_reg} using Algorithm~\ref{alg:pod-based} or \ref{alg:am-based}
\STATE Approximate $\dfrac{\text{d}}{\text{d}t}{\hat{\mathbf{S}}}$ from $\hat{\mathbf{S}}$ \hfill \Comment{Time derivative approximation}
\STATE Choose $\lambda_1$, $\lambda_2$, $\lambda_3$ \hfill \Comment{Hyperparameter optimization}
\STATE Solve OpInf regression problem \eqref{eq:opinf_nonlinear} to find $\hat{\mathbf{c}}$, $\hat{\mathbf{A}}$, $\hat{\mathbf{H}}$, $\hat{\mathbf{P}}$.
\end{algorithmic}
\end{algorithm}
\vspace{-1em}
\end{figure}

\section{Numerical experiments}
\label{sec:numerical_experiments}

In this section, we discuss application of the OpInf model reduction methods described in Section~\ref{sec:roms} to several dynamical systems. 
We compare the OpInf approach from Peherstorfer and Willcox\cite{PEHERSTORFER2016196} (Algorithm~\ref{alg:standard_opinf}) and the nonlinear-manifold-based OpInf approaches proposed above: MPOD-OpInf and MAM-OpInf (see Algorithm~\ref{alg:proposed_opinf}). 
We report numerical results for benchmark problems involving the  Allen-Cahn equation, the Korteweg-de Vries equation, and a cylinder flow problem. 

\subsection{Practical considerations}
\label{subsec:practical_considerations}

The dimensionality of a reduced-order model is typically informed by the representation error of the training data used in its construction. 
For nonlinear approximations of the form \eqref{eq:nonlinear_approx} we compute the metric\cite{GEELEN2023115717}
\begin{equation}
\epsilon_r = \dfrac{\| \mathbf{V}\hat{\mathbf{S}} + \overline{\mathbf{V}}\boldsymbol{\Xi}\mathbf{g}(\hat{\mathbf{S}}) \|_F^2}{ \| \mathbf{S}-\mathbf{S}_\text{ref} \|_F^2 }.
\label{eq:rit}
\end{equation}
When $\overline{\mathbf{V}}$ is zero (as in the linear-subspace OpInf approach) we recover the well-known expression 
\begin{equation}
\epsilon_r = \sum_{j=1}^r \sigma_j^2 / \sum_{j=1}^k \sigma_j^2,
\label{eq:snapshot_energy}
\end{equation} 
where $\sigma_j$ denotes the $j$th singular value of the mean-centered snapshot matrix $\mathbf{S}-\mathbf{S}_\text{ref}$. This indicator is commonly referred to as the (cumulative) snapshot energy captured by the basis. The dimensionality of the reduced-order models, $r$, and the number of orthogonal basis vectors in the nonlinear part of the approximation, $q$, are user-specified. In the following these values are chosen based on the singular value decay of the shifted snapshot matrices. The polynomial order, $p$, of approximation \eqref{eq:nonlinear_approx} is chosen, in accordance with $r$ and $q$, based on the ability to represent data in the reduced-state coordinate system to sufficient accuracy.

The primary error metric used in numerical experiments is the relative error in the states, namely, $\| \mathbf{S}-\boldsymbol{\Gamma}(\hat{\mathbf{S}}) \|_F/ \| \mathbf{S}-\mathbf{S}_\text{ref} \|_F$. 
Initial guesses for $\mathbf{V}, \overline{\mathbf{V}}, \boldsymbol{\Xi}, \hat{\mathbf{S}}$ in the alternating minimization method from Algorithm~\ref{alg:am-based} are obtained from the POD-based representation learning method of Algorithm~\ref{alg:pod-based}. 
The iterative process from Algorithm~\ref{alg:am-based} is terminated when the change in the relative snapshot energy $\epsilon_r$ \eqref{eq:rit} on consecutive iterations falls below $10^{-3}$. 
Note that this termination criterion for the alternating minimization scheme is not indicative of the snapshot energy captured by a given number of basis vectors and reduced-state coefficients. One might change this threshold based on model size or desired optimization tolerance.
The function tolerance for the nonlinear least-squares solver used in solving \eqref{eq:opt_low_dim_coord_system} is set to $10^{-9}$. 
It is noted that the MAM-OpInf approach calls for a representation of the initial condition in the reduced-state
coordinate system. This involves computing $r$ coefficients through a nonlinear least-squares problem which is carried out at the start of each reduced model evaluation.

The regularization parameters $\lambda_i$ with $i=1,2,3$ in \eqref{eq:opinf_nonlinear} and the regularization parameter $\gamma$ in representation learning problem \eqref{eq:optim_problem_reg} are calibrated through a grid search conducted over a predetermined range of candidate values. We seek the parameter combination that minimizes the relative state error over the available training data.\cite{doi:10.1080/03036758.2020.1863237} This requires that the reduced-order model be evaluated for every parameter combination. We also note that $\gamma$ does not (explicitly) show up in OpInf problem \eqref{eq:opinf_nonlinear} as it appears only upon computing a representation of the data in the original state-space through \eqref{eq:nonlinear_approx}.
The time derivatives of the projected snapshot data in the OpInf regression problems are estimated via a fourth-order finite difference approximation.

\subsection{The Allen-Cahn equation}
\label{subsec:allen-cahn}

While the Allen-Cahn model was originally conceived to describe the motion of anti-phase boundaries in metallic alloys,\cite{ALLEN19791085} it has become prototypical for describing phase separation and interfacial dynamics in many application domains. The equation has also been studied in the context of model reduction.\cite{SONG2016213, a13060148}
We consider the Allen-Cahn equation 
\begin{equation}
    \partial_t s = \kappa \partial_{x}^2 s + s - s^3
    \label{eq:allen_cahn}
\end{equation}
in the domain $x \in [-1,1]$ with Dirichlet boundary conditions $s(-1,t)=-1; s(1,t)=1$ and initial condition
\begin{equation}
    s(x,0) = \mu x + (1-\mu)\sin(-1.5 \pi x),
\end{equation}
in which the parameter $\mu$ varies uniformly on the range $[0.5, 0.6]$. The operators $\partial_x$ and $\partial_t$ in \eqref{eq:allen_cahn} denote partial differentiation with respect to space and time, respectively.
The interface parameter $\kappa \in \mathbb{R}^+$ is positive constant which represents the thickness of the interface that separates the two phases.

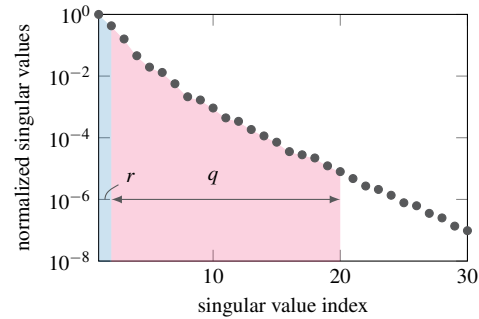
\begin{figure}[!tbp]
\centering \footnotesize
\begin{tikzpicture}
\begin{semilogyaxis}[width=.75\linewidth, height=.5625\linewidth, xlabel={singular value index}, ylabel={normalized singular values}, ytickten={0,-2,...,-8}, xmin=1, xmax=30, ymin=1e-8, ymax=1e0, legend style={draw=none}, legend cell align={left}, mark=none, grid=none]
\addplot[only marks, Black!80, mark=*, mark size=1.5pt] table[x=index_a, y=norm_sv_a, col sep=comma] {./data/sv_1d_ac_chaos.csv};
\addplot+[Black!0, mark=none, name path=A] table[x=index, y=norm_sv, col sep=comma, mark=none] {./data/sv_1d_ac_chaos.csv};
\addplot+[Black!0, mark=none, name path=B] table[x=index, y=zero, col sep=comma, mark=none] {./data/sv_1d_ac_chaos.csv};
\addplot+[Black!0, mark=none, name path=C] table[x=index2, y=norm_sv2, col sep=comma, mark=none] {./data/sv_1d_ac_chaos.csv};
\addplot+[Black!0, mark=none, name path=D] table[x=index2, y=zero2, col sep=comma, mark=none] {./data/sv_1d_ac_chaos.csv};
\addplot+[NavyBlue, fill opacity=0.20] fill between[of=A and B]; % filling
\addplot+[WildStrawberry, fill opacity=0.20] fill between[of=C and D]; % filling
\node[draw=none] (box) at (3.5,5e-6) {$r$};
\coordinate (a) at (1.5,1e-6) ;
\draw [Black!80](a) to[out=90, in=-150] (box) ;
\node[draw=none] at (10,5e-6) {$q$};
\coordinate (c) at (2,1e-6) ;
\coordinate (d) at (20,1e-6) ;
\draw [latex-latex, Black!80](c) -- (d) ;
\end{semilogyaxis}
\end{tikzpicture}
\caption{Normalized singular values of the centered snapshot matrix for the Allen-Cahn problem. The blue and red areas denote the singular values whose corresponding left singular vectors are columns in $\mathbf{V}$ and $\overline{\mathbf{V}}$, respectively, in the MPOD-OpInf formulation.}
\label{fig:svs_ac}
\end{figure}

\begin{figure}[!tbp]
    \centering \footnotesize
\begin{subfigure}{.9\linewidth}
    \begin{overpic}[width=\linewidth]{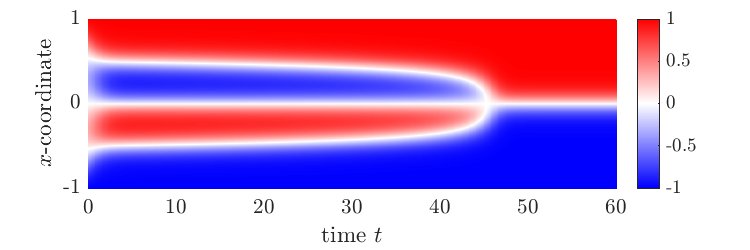}
    \end{overpic}
    \caption{Reference}
\end{subfigure}    
\begin{subfigure}{.9\linewidth}
    \begin{overpic}[width=\linewidth]{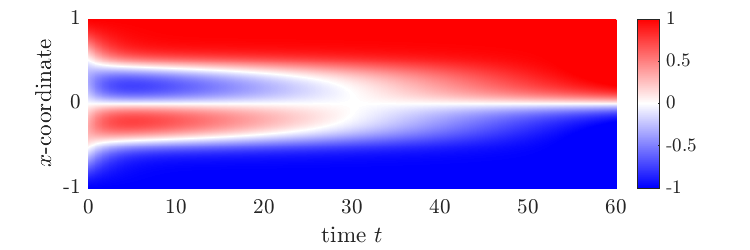}
    \end{overpic}
    \caption{Linear-subspace OpInf (POD): $r=2$}
\end{subfigure}    
\begin{subfigure}{.9\linewidth}
    \begin{overpic}[width=\linewidth]{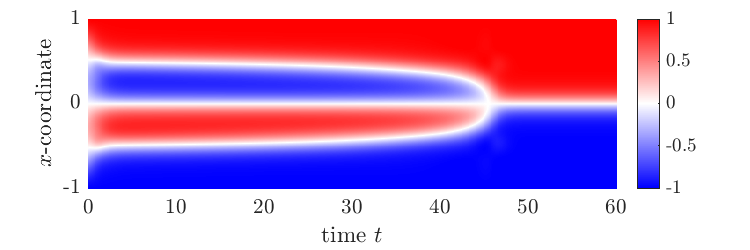}
    \end{overpic}
    \caption{MPOD-OpInf: $r=2$, $p=4$}
\end{subfigure}    
\caption{Comparison of the reference solution at parameter value $\mu=0.5127$ (top) with the reconstructions from the two-equation OpInf (middle) and MPOD-OpInf (bottom) models for the Allen-Cahn model.}
\label{fig:allen_cahn_reconstruction}
\end{figure}

\begin{figure}[!tbp]
    \centering \footnotesize
\begin{subfigure}{.9\linewidth}
    \begin{overpic}[width=\linewidth]{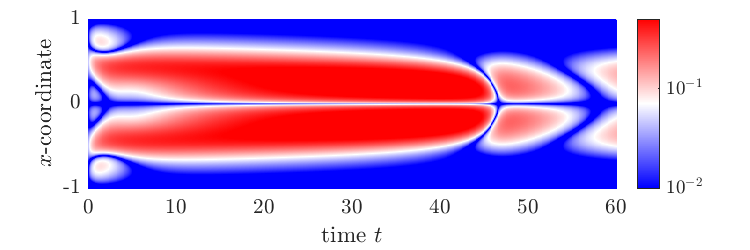}
    \end{overpic}
    \caption{Linear-subspace OpInf (POD): $r=2$}
\end{subfigure}    
\begin{subfigure}{.9\linewidth}
    \begin{overpic}[width=\linewidth]{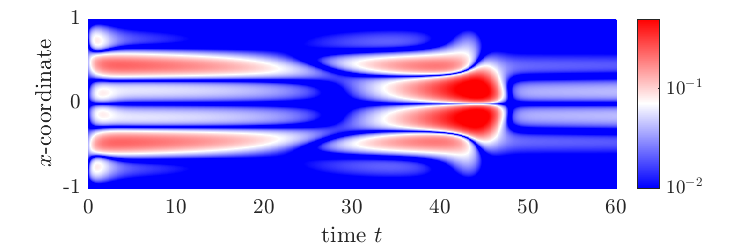}
    \end{overpic}
    \caption{MPOD-OpInf: $r=2$, $p=2$}
\end{subfigure}    
\begin{subfigure}{.9\linewidth}
    \begin{overpic}[width=\linewidth]{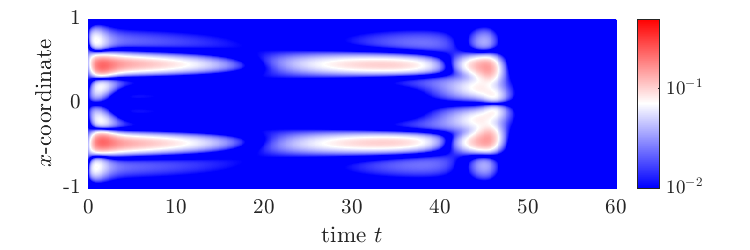}
    \end{overpic}
    \caption{MPOD-OpInf: $r=2$, $p=3$}
\end{subfigure}    
\begin{subfigure}{.9\linewidth}
    \begin{overpic}[width=\linewidth]{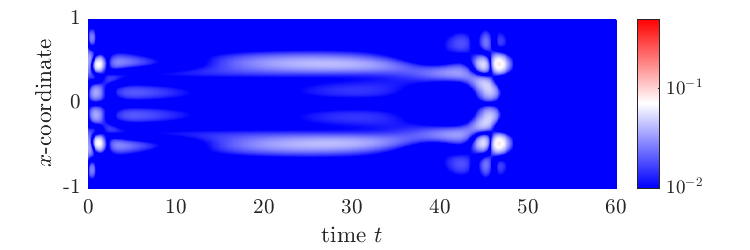}
    \end{overpic}
    \caption{MPOD-OpInf: $r=2$, $p=4$}
\end{subfigure}    
\caption{Pointwise error in the reconstructions for the test trajectory at parameter value $\mu=0.5127$ for the Allen-Cahn model.}
\label{fig:allen_cahn_error}
\end{figure}

While the Allen-Cahn equation \eqref{eq:allen_cahn} is characterized by linear and cubic terms, adding an auxiliary variable yields system dynamics with the desired quadratic model structure \eqref{eq:fom}. We apply the transformation\cite{QIAN2020132401, doi:10.1137/14097255X}
\begin{equation}
    \mathcal{T} : 
    \begin{pmatrix}
        s \\
    \end{pmatrix}
    \mapsto 
    \begin{pmatrix}
        s \\
        s^2
    \end{pmatrix}
    \equiv
    \begin{pmatrix}
        w_1 \\
        w_2
    \end{pmatrix}  
\end{equation}
The lifted system is then given by
\begin{equation}
\begin{aligned}
    \partial_t w_1 &= \epsilon \partial_{x}^2 w_1 + w_1 + w_1w_2 \\
    \partial_t w_2 &= 2w_1 \partial_t w_1 \\
                   &= 2w_1 \left( \epsilon \partial_{x}^2 w_1 + w_1 + w_1w_2 \right) \\
                   &= 2 \epsilon w_1 \partial_{x}^2 w_1 + 2w_2 + 2 (w_2)^2,
    \label{eq:allen_cahn_lifted}
\end{aligned}
\end{equation}
which contains only quadratic nonlinear dependencies on the state. It is also important to note that no approximations are invoked in the process of lifting \eqref{eq:allen_cahn} to \eqref{eq:allen_cahn_lifted}. 

State data are computed on a uniform spatial grid consisting of $n = 512$ grid points. 
The state snapshot data are generated from three simulations of the Allen-Cahn model, corresponding to the parameters $\mu = [0.50, 0.55, 0.60]$. For testing, ten more trajectories are generated with parameter $\mu$ drawn uniformly at random from the interval $[0.5,0.6]$.
The data are recorded every $0.1$ time units up to time $T=60$, yielding 600 snapshots per trajectory, for a total of 1,800. 
The lifted snapshot matrix is centered by the mean initial condition across the test parameters. 

Fig.~\ref{fig:svs_ac} shows the decay of the normalized singular values, where normalized means that the first normalized singular value equals 1. The first two POD modes contain 97.7\% of the energy in the lifted state data, but a total of $r+q=20$ POD modes are needed to drive the projection error in the training data below $10^{-5}$. 
We consider two-equation reduced-order models constructed from OpInf methodologies. 
For the MPOD-OpInf formulation, the basis matrix $\mathbf{V}$ contains the first $r=2$ POD modes, with the remaining $q=20-r=18$ POD modes captured in $\overline{\mathbf{V}}$.
The regularization parameter in the representation learning problem in \eqref{eq:optim_problem_reg} is chosen to be $\gamma = 10^{-2}$. 
Value for the regularization parameters  $\lambda_i$, $i=1,2,3$ in the OpInf problems are found by minimizing the relative state error across all the training parameters $\mu$. 

\begin{table}[!tbp]
    \centering
    \caption{Median of relative state error \eqref{eq:rit} in the training and test problems across the parameters for the Allen-Cahn problem.}
    \label{tab:allen_cahn_error}
    \renewcommand{\arraystretch}{1.25}%
    \begin{tabular}{c|c|c|}
    \cline{2-3}
         & \bfseries Training & \bfseries Testing \\
    \hline
        \multicolumn{1}{|c|}{\bfseries Linear-subspace OpInf}        & $3.599 \times 10^{-1}$ & $3.424\times 10^{-1}$ \\
        \multicolumn{1}{|c|}{\bfseries MPOD-OpInf ($p=2$)} & $1.823 \times 10^{-1}$ & $1.552 \times 10^{-1}$ \\
        \multicolumn{1}{|c|}{\bfseries MPOD-OpInf ($p=3$)} & $5.091 \times 10^{-2}$ & $4.368 \times 10^{-2}$ \\
        \multicolumn{1}{|c|}{\bfseries MPOD-OpInf ($p=4$)} & $2.567 \times 10^{-2}$ & $2.552 \times 10^{-2}$ \\
    \hline
    \end{tabular}
\end{table}

Fig.~\ref{fig:allen_cahn_reconstruction} shows the reconstructed trajectories at a random test parameter $\mu=0.5127$ for the linear-subspace OpInf formulation and its manifold-based counterpart using fourth-order polynomial embeddings. 
The MPOD-OpInf model represents the phase separation process over time more accurately, both in the material phases and the interface dynamics. 
Pointwise errors in the reconstructions of the original state data for the nonlinear manifold models are shown in Fig.~\ref{fig:allen_cahn_error}. 
We note from these figures that increasing the degree of the polynomial embeddings can improve the  predictive capabilities. 
(Performance of a linear-subspace OpInf reduced-order model is shown for reference.) 
The relative state errors for the reduced-order models are tabulated in Table \ref{tab:allen_cahn_error}. 
It can be seen that the MPOD-OpInf formulation outperforms the linear-subspace OpInf formulation in the training regime as well as in the predictive setting.

\subsection{The Korteweg-de Vries equation}
\label{subsec:kdv}

The second numerical experiment is concerned with traveling wave physics.\cite{mendible2020dimensionality} We consider a single propagating soliton in a one-dimensional domain with periodic boundary conditions. The evolution of the wave field $s$ in the space-time domain $[-\pi,\pi]\times[0,T]$ is obtained from the Korteweg-de Vries equation
\begin{equation}
    \partial_t s = -\alpha s \partial_x s - \beta \partial_{x}^3 s.
\end{equation}
The initial condition is given by $s_0(x) = 1 + 24 \sech^2 \left( \sqrt{8}x \right)$. 
We use an equidistant computational grid consisting of 256 evenly spaced points in space. 
State data are saved every $0.0002$ time units. 
We choose a final time $T=1$ and model constants $\alpha=4$ and $\beta=1$.

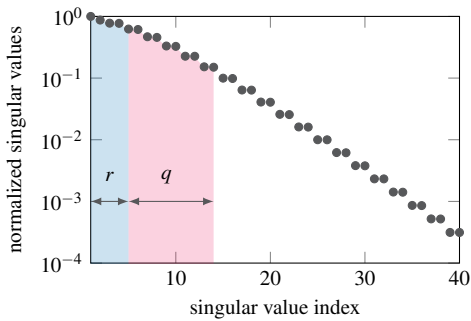
\begin{figure}[!tbp]
\centering \footnotesize
\begin{tikzpicture}
\begin{semilogyaxis}[width=.75\linewidth, height=.5625\linewidth, xlabel={singular value index}, ylabel={normalized singular values}, ytickten={0,-1,...,-4}, xmin=1, xmax=40, ymin=1e-4, ymax=1e0, legend style={draw=none}, legend cell align={left}, mark=none, grid=none]
\addplot[only marks, Black!80, mark=*, mark size=1.5pt] table[x=index_a, y=norm_sv_a, col sep=comma] {./data/sv_1d_chaos.csv};
\addplot+[Black!0, mark=none, name path=A] table[x=index, y=norm_sv, col sep=comma, mark=none] {./data/sv_1d_chaos.csv};
\addplot+[Black!0, mark=none, name path=B] table[x=index, y=zero, col sep=comma, mark=none] {./data/sv_1d_chaos.csv};
\addplot+[Black!0, mark=none, name path=C] table[x=index2, y=norm_sv2, col sep=comma, mark=none] {./data/sv_1d_chaos.csv};
\addplot+[Black!0, mark=none, name path=D] table[x=index2, y=zero2, col sep=comma, mark=none] {./data/sv_1d_chaos.csv};
\addplot+[NavyBlue, fill opacity=0.20] fill between[of=A and B]; % filling
\addplot+[WildStrawberry, fill opacity=0.20] fill between[of=C and D]; % filling
\node[draw=none] at (3,2.5e-3) {$r$};
\coordinate (a) at (1,1e-3) ;
\coordinate (b) at (5,1e-3) ;
\draw [latex-latex, Black!80](a) -- (b) ;
\node[draw=none] at (9,2.5e-3) {$q$};
\coordinate (c) at (5,1e-3) ;
\coordinate (d) at (14,1e-3) ;
\draw [latex-latex, Black!80](c) -- (d) ;
\end{semilogyaxis}
\end{tikzpicture}
\caption{Normalized singular values of the mean-subtracted snapshot matrix for the Korteweg-de Vries problem. The blue and red areas denote the singular values whose corresponding left singular vectors are columns in $\mathbf{V}$ and $\overline{\mathbf{V}}$, respectively, in the MPOD-OpInf formulation.}
\label{fig:svs_kdv}
\end{figure}

\begin{figure}[!tbp]
    \centering
\begin{subfigure}{.9\linewidth}
    \includegraphics[width=\linewidth]{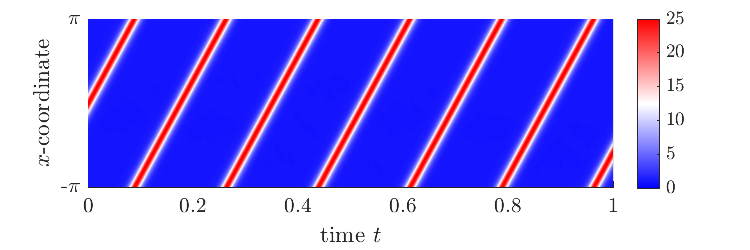}
    \caption{Reference}
\end{subfigure}    
\begin{subfigure}{.9\linewidth}
    \includegraphics[width=\linewidth]{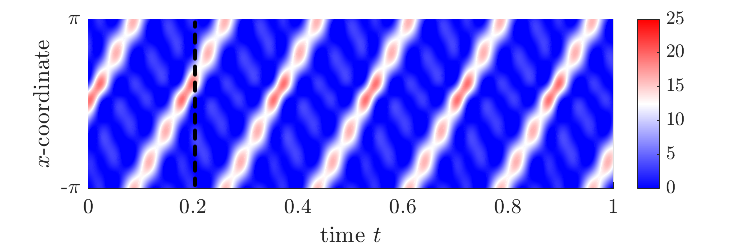}
    \caption{Linear-subspace OpInf, $r=5$}
\end{subfigure}    
\begin{subfigure}{.9\linewidth}
    \includegraphics[width=\linewidth]{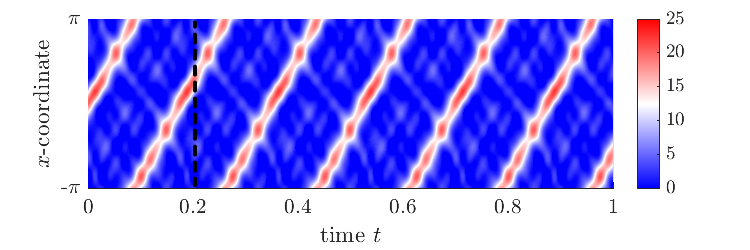}
    \caption{MPOD-OpInf, $r=5$, $p=2$}
\end{subfigure}    
\begin{subfigure}{.9\linewidth}
    \includegraphics[width=\linewidth]{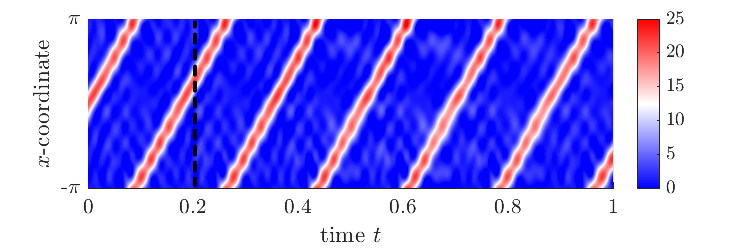}
    \caption{MAM-OpInf, $r=5$, $p=2$}
\end{subfigure}    
\caption{Plots of the reference solution and OpInf-produced predictions for the Korteweg-de Vries equation over the time window $t \in [0,T]$ at a reduced basis dimension of $r=5$. The end of the training window is indicated by the dashed line.}
\label{fig:kdv_r12}
\end{figure}

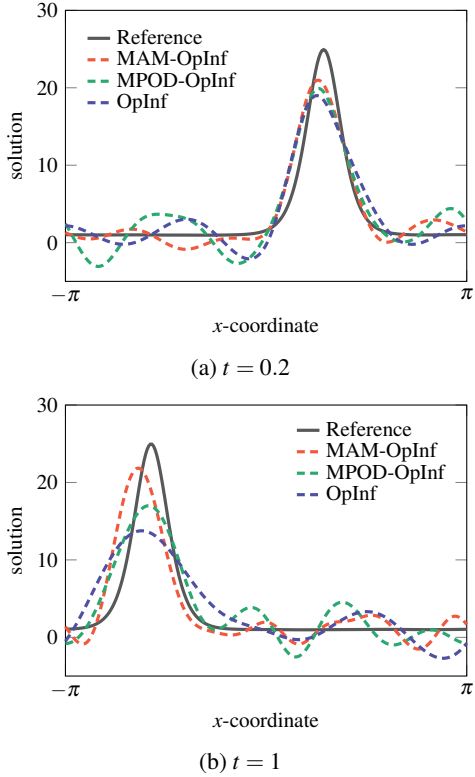
\begin{figure}[!tbp]
\centering \footnotesize
\begin{subfigure}[t]{.8\linewidth}
\begin{tikzpicture}
\begin{axis}[width=\linewidth, height=.75\linewidth, xlabel={$x$-coordinate}, ylabel={solution}, xtick={-3.1415,3.1415}, xticklabels={$-\pi$,$\pi$}, xmin=-3.1416, xmax=3.1416, grid=none, ymin=-5, ymax=30, xlabel near ticks, ylabel near ticks, legend style={ fill=none, draw=none}, legend pos=north west, legend cell align={left}, legend style={row sep=-2.5pt}]
\addplot+[solid, very thick, Black!80, mark=none] table[x=x, y=T20, col sep=comma, mark=none] {./data/kdv_snapshots2.csv};
\addlegendentry{Reference}
\addplot+[densely dashed, very thick, Red!80, mark=none] table[x=x, y=MAM_T20, col sep=comma, mark=none] {./data/kdv_snapshots2.csv};
\addlegendentry{MAM-OpInf}
\addplot+[densely dashed, very thick, ForestGreen!80, mark=none] table[x=x, y=MPOD_T20, col sep=comma, mark=none] {./data/kdv_snapshots2.csv};
\addlegendentry{MPOD-OpInf}
\addplot+[densely dashed, very thick, Blue!80, mark=none] table[x=x, y=POD_T20, col sep=comma, mark=none] {./data/kdv_snapshots2.csv};
\addlegendentry{OpInf}
\end{axis}
\end{tikzpicture}
\caption{$t=0.2$}
\end{subfigure}
\begin{subfigure}[t]{.8\linewidth}
\begin{tikzpicture}
\begin{axis}[width=\linewidth, height=.75\linewidth, xlabel={$x$-coordinate}, ylabel={solution}, xtick={-3.1415,3.1415}, xticklabels={$-\pi$,$\pi$}, xmin=-3.1416, xmax=3.1416, grid=none, ymin=-5, ymax=30, xlabel near ticks, ylabel near ticks, legend style={ fill=none, draw=none}, legend cell align={left}, legend style={row sep=-2.5pt}, clip=false]
\addplot+[solid, very thick, Black!80, mark=none] table[x=x, y=T100, col sep=comma, mark=none] {./data/kdv_snapshots2.csv};
\addlegendentry{Reference}
\addplot+[densely dashed, very thick, Red!80, mark=none] table[x=x, y=MAM_T100, col sep=comma, mark=none] {./data/kdv_snapshots2.csv};
\addlegendentry{MAM-OpInf}
\addplot+[densely dashed, very thick, ForestGreen!80, mark=none] table[x=x, y=MPOD_T100, col sep=comma, mark=none] {./data/kdv_snapshots2.csv};
\addlegendentry{MPOD-OpInf}
\addplot+[densely dashed, very thick, Blue!80, mark=none] table[x=x, y=POD_T100, col sep=comma, mark=none] {./data/kdv_snapshots2.csv};
\addlegendentry{OpInf}
\end{axis}
\end{tikzpicture}
\caption{$t=1$}
\end{subfigure}
\caption{Reference and OpInf solution snapshots from the Korteweg-de Vries experiment at the end of the training regime ($t=0.2$) and predictions for the final time $T$. The plots are given for models of dimension $r=5$.}
\label{fig:kdv_snapshots}
\end{figure}

To learn the nonlinear manifolds and train our data-driven reduced-order models, 1001 snapshots of the solution are collected uniformly across the time interval $t \in [0, 0.2]$. 
The snapshot matrix under consideration is centered by its column-averaged (thus time-averaged) mean value; the decay of its singular values is shown in Fig.~\ref{fig:svs_kdv}. 
The dynamics of the system can be captured well with only 14
modes capturing 99.3\% of the cumulative snapshot energy.
The values of $r$ (the number of columns of $\mathbf{V}$) and $q$ (the number of columns of $\overline{\mathbf{V}}$) are then chosen so that $r+q=14$. 
The regularization parameter in the representation learning problem in \eqref{eq:optim_problem_reg} is chosen so that $\gamma = 10^{-3}$. 
We now consider the performance of the reduced-order models in both the training regime as well as in the predictive setting.

\begin{figure}[!tbp]
\centering \scriptsize
\begin{subfigure}[t]{.495\linewidth}
\begin{tikzpicture}
\begin{axis}[width=1.05\linewidth, xlabel={$x$-coordinate}, ylabel={mode shape}, xtick={-3.1415,3.1415}, xticklabels={$-\pi$,$\pi$}, xmin=-3.1416, xmax=3.1416, grid=none, ytickten={-0.1,0.1}, ymin=-0.17, ymax=0.17, xlabel near ticks, ylabel near ticks, legend style={ fill=none, draw=none}, legend cell align={left}, legend style={row sep=-2.5pt}, clip=false]
\addplot+[solid, very thick, Black!80, mark=none] table[x=x, y=POD1, col sep=comma, mark=none] {./data/kdv_procrustes_modes2.csv};
\addplot+[densely dashed, very thick, Red!80, mark=none] table[x=x, y=Pro1, col sep=comma, mark=none] {./data/kdv_procrustes_modes2.csv};
\end{axis}
\end{tikzpicture}
\caption{Basis vector 1}
\end{subfigure}
\begin{subfigure}[t]{.495\linewidth}
\begin{tikzpicture}
\begin{axis}[width=1.05\linewidth, xlabel={$x$-coordinate}, ylabel={mode shape}, xtick={-3.1415,3.1415}, xticklabels={$-\pi$,$\pi$}, xmin=-3.1416, xmax=3.1416, grid=none, ytickten={-0.1,0.1}, ymin=-0.17, ymax=0.17, xlabel near ticks, ylabel near ticks, legend style={ fill=none, draw=none}, legend cell align={left}, legend style={row sep=-2.5pt}, clip=false]
\addplot+[solid, very thick, Black!80, mark=none] table[x=x, y=POD2, col sep=comma, mark=none] {./data/kdv_procrustes_modes2.csv};
\addplot+[densely dashed, very thick, Red!80, mark=none] table[x=x, y=Pro2, col sep=comma, mark=none] {./data/kdv_procrustes_modes2.csv};
\end{axis}
\end{tikzpicture}
\caption{Basis vector 2}
\end{subfigure} \\[1em]
\begin{subfigure}[t]{.495\linewidth}
\begin{tikzpicture}
\begin{axis}[width=1.05\linewidth, xlabel={$x$-coordinate}, ylabel={mode shape}, xtick={-3.1415,3.1415}, xticklabels={$-\pi$,$\pi$}, xmin=-3.1416, xmax=3.1416, grid=none, ytickten={-0.1,0.1}, ymin=-0.17, ymax=0.17, xlabel near ticks, ylabel near ticks, legend style={ fill=none, draw=none}, legend cell align={left}, legend style={row sep=-2.5pt}, clip=false]
\addplot+[solid, very thick, Black!80, mark=none] table[x=x, y=POD3, col sep=comma, mark=none] {./data/kdv_procrustes_modes2.csv};
\addplot+[densely dashed, very thick, Red!80, mark=none] table[x=x, y=Pro3, col sep=comma, mark=none] {./data/kdv_procrustes_modes2.csv};
\end{axis}
\end{tikzpicture}
\caption{Basis vector 3}
\end{subfigure} 
\begin{subfigure}[t]{.495\linewidth}
\begin{tikzpicture}
\begin{axis}[width=1.05\linewidth, xlabel={$x$-coordinate}, ylabel={mode shape}, xtick={-3.1415,3.1415}, xticklabels={$-\pi$,$\pi$}, xmin=-3.1416, xmax=3.1416, grid=none, ytickten={-0.1,0.1}, ymin=-0.17, ymax=0.17, xlabel near ticks, ylabel near ticks, legend style={ fill=none, draw=none}, legend cell align={left}, legend style={row sep=-2.5pt}, clip=false]
\addplot+[solid, very thick, Black!80, mark=none] table[x=x, y=POD4, col sep=comma, mark=none] {./data/kdv_procrustes_modes2.csv};
\addplot+[densely dashed, very thick, Red!80, mark=none] table[x=x, y=Pro4, col sep=comma, mark=none] {./data/kdv_procrustes_modes2.csv};
\end{axis}
\end{tikzpicture}
\caption{Basis vector 4}
\end{subfigure}
\caption{First four basis vectors computed using POD (black curves) and those obtained by means of orthogonal Procrustes problem \eqref{eq:procrustes} in Algorithm~\ref{alg:am-based} (red curves) for the Korteweg-de Vries problem. For the latter curves, the modes are computed for models of size $r=5$ with quadratic embeddings ($p=2$).}
\label{fig:kdv_basis_vectors}
\end{figure}
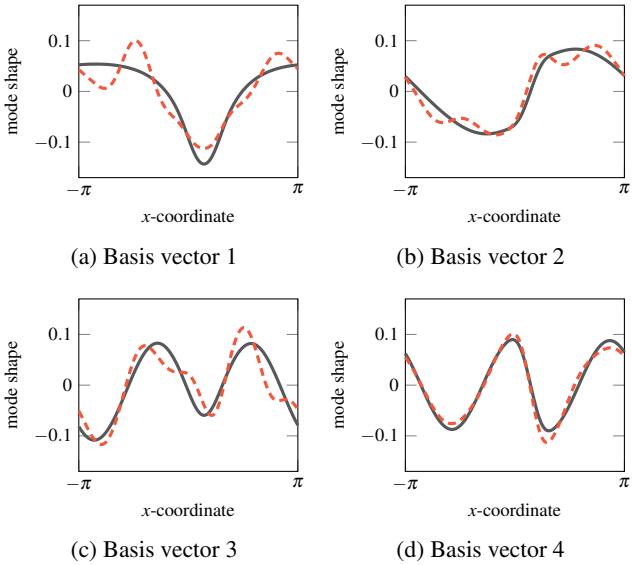

Visual comparisons of the reference solution and the solutions produced by the OpInf models are shown in Fig.~\ref{fig:kdv_r12}. 
The MPOD-OpInf and MAM-OpInf models use polynomial embeddings of degree $p=2$. 
Even at $r = 5$ basis vectors, the space-time evolution of the propagating soliton is relatively well captured. 
The OpInf, MPOD-OpInf, and MAM-OpInf models account for 73.6\%, 82.0\% and 93.1\%, respectively, of the cumulative energy \eqref{eq:rit} in the system. 
The relative state errors for these formulations, over the window of training, equals 51.4\%, 42.5\% and 30.0\%. 
Predictably, the addition of quadratic terms in the state approximation yields an increase in accuracy. 
By accounting for orthogonal transformations of the basis, as induced by the alternating minimization procedure, the contrast in accuracy becomes more pronounced, as the MAM-OpInf model is more accurate than the MPOD-OpInf variant (despite these two models having the exact same online computational expense).

Outside the range of their training data (that is $t\in[0.2,1]$) the relative state errors are 52.9\%, 43.3\%, and 35.1\%, for the OpInf, MPOD-OpInf, and MAM-OpInf models respectively. This experiment demonstrates the potential of both the MPOD-OpInf and MAM-OpInf methods to outperform linear-subspace OpInf methods in the training and (especially) in the predictive regime. These results are corroborated by comparing the solution snapshots at the end of training ($t=0.2$) and at final time ($t=1$) for the different methods (see Fig.~\ref{fig:kdv_snapshots}). While there is some error associated with the prediction of the soliton's exact spatial location, the MPOD-OpInf, and MAM-OpInf formulation are better suited for capturing the soliton's representation over time.

The first four basis vectors (corresponding to the dominant left singular vectors) for state approximations of dimension $r=5$ of order $p=2$ are shown in Fig.~\ref{fig:kdv_basis_vectors}. 
While both the MPOD-OpInf and MAM-OpInf approaches produce an orthogonal set of vectors, the alternating minimization approach can be seen to incorporate some of the small-scale solution features into the dominant modes. 

\begin{figure}[!tbp]
    \centering
\begin{subfigure}{.9\linewidth}
    \includegraphics[width=\linewidth]{./figures/REF.png}
    \caption{Reference}
\end{subfigure}    
\begin{subfigure}{.9\linewidth}
    \includegraphics[width=\linewidth]{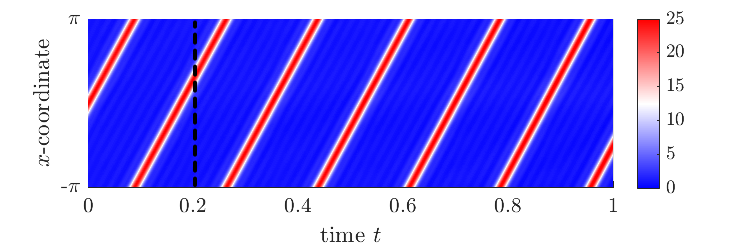}
    \caption{Linear-subspace OpInf, $r=16$}
\end{subfigure}    
\begin{subfigure}{.9\linewidth}
    \includegraphics[width=\linewidth]{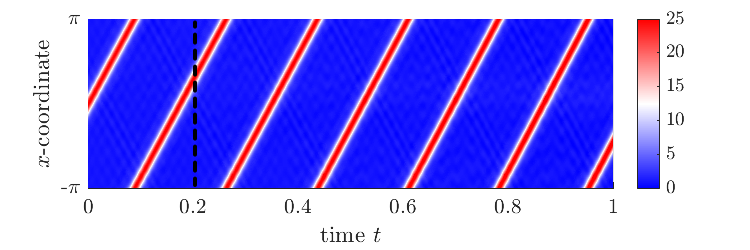}
    \caption{MPOD-OpInf, $r=16$, $p=2$}
\end{subfigure}    
\begin{subfigure}{.9\linewidth}
    \includegraphics[width=\linewidth]{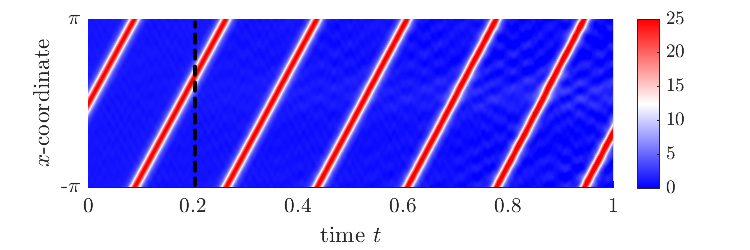}
    \caption{MAM-OpInf, $r=16$, $p=2$}
\end{subfigure}    
\caption{Plots of the reference solution and OpInf-produced predictions for the Korteweg-de Vries equation over the time window $t \in [0,T]$ at a reduced basis dimension of $r=16$. The end of the training window is indicated by the dashed line.}
\label{fig:kdv_r16}
\end{figure}

\begin{figure}[!tbp]
\centering \footnotesize
\begin{subfigure}[t]{.8\linewidth}
\begin{tikzpicture}
\begin{axis}[width=\linewidth, height=.75\linewidth, xlabel={$x$-coordinate}, ylabel={solution}, xtick={-3.1415,3.1415}, xticklabels={$-\pi$,$\pi$}, xmin=-3.1416, xmax=3.1416, grid=none, ymin=-5, ymax=30, xlabel near ticks, ylabel near ticks, legend style={ fill=none, draw=none}, legend cell align={left}, legend style={row sep=-2.5pt}, clip=false]
\addplot+[solid, very thick, Black!80, mark=none] table[x=x, y=T100, col sep=comma, mark=none] {./data/kdv_snapshots3.csv};
\addlegendentry{Reference}
\addplot+[densely dashed, very thick, Blue!80, mark=none] table[x=x, y=POD_T100, col sep=comma, mark=none] {./data/kdv_snapshots3.csv};
\addlegendentry{OpInf}
\addplot+[densely dashed, very thick, ForestGreen!80, mark=none] table[x=x, y=MPOD_T100, col sep=comma, mark=none] {./data/kdv_snapshots3.csv};
\addlegendentry{MPOD-OpInf}
\addplot+[densely dashed, very thick, Red!80, mark=none] table[x=x, y=MAM_T100, col sep=comma, mark=none] {./data/kdv_snapshots3.csv};
\addlegendentry{MAM-OpInf}
\end{axis}
\end{tikzpicture}
\end{subfigure}
\caption{Reference and OpInf solution snapshots from the Korteweg-de Vries experiment at final time $t=1$. The plots are given for models of dimension $r=16$.}
\label{fig:kdv_snapshots3}
\end{figure}
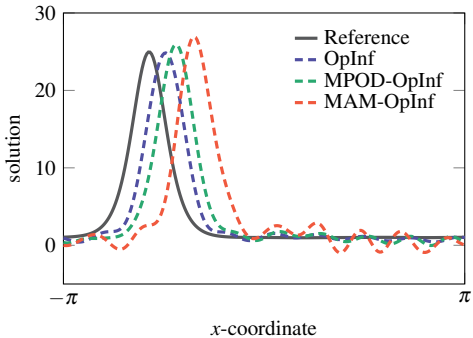

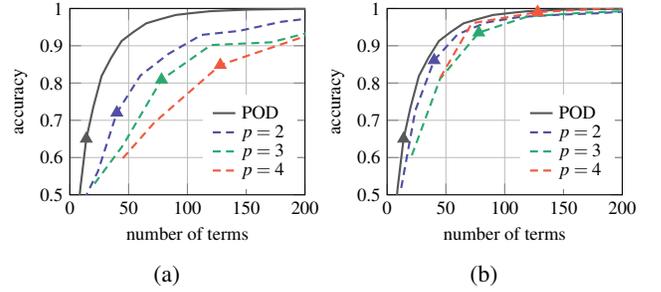
\begin{figure}[!tbp]
\centering \scriptsize
\begin{subfigure}[t]{.495\linewidth}
\begin{tikzpicture}
\begin{axis}[width=1.1\linewidth, xlabel={number of terms}, ylabel={accuracy}, xmin=0, xmax=200, ytick={0.5,0.6,0.7,0.8,0.9,1}, yticklabels={0.5,0.6,0.7,0.8,0.9,1}, grid=both, ymin=0.5, ymax=1.0, xlabel near ticks, ylabel near ticks, legend style={draw=none}, legend pos=south east, legend cell align={left}, legend style={row sep=-2.5pt}]
\addplot+[solid, thick, Black!80, mark=none] table[x=cost1, y=POD, col sep=comma, mark=none] {./data/kdv_cost.csv};
\addlegendentry{POD}
\addplot+[densely dashed, thick, Blue!80, mark=none] table[x=cost2, y=MPOD2, col sep=comma, mark=none] {./data/kdv_cost.csv};
\addlegendentry{$p=2$}
\addplot+[densely dashed, thick, ForestGreen!80, mark=none] table[x=cost3, y=MPOD3, col sep=comma, mark=none] {./data/kdv_cost.csv};
\addlegendentry{$p=3$}
\addplot+[densely dashed, thick, Red!80, mark=none] table[x=cost4, y=MPOD4, col sep=comma, mark=none] {./data/kdv_cost.csv};
\addlegendentry{$p=4$}
\node[Black!80, mark size=2.5pt] (A) at (14,0.650245) {\pgfuseplotmark{triangle*}};
\node[Blue!80, mark size=2.5pt] (B) at (40,0.720444) {\pgfuseplotmark{triangle*}};
\node[ForestGreen!80, mark size=2.5pt] (C) at (78,0.808955) {\pgfuseplotmark{triangle*}};
\node[Red!80, mark size=2.5pt] (D) at (128,0.848390) {\pgfuseplotmark{triangle*}};
\end{axis}
\end{tikzpicture}
\caption{}
\end{subfigure}
\hspace{-0.75em}
\begin{subfigure}[t]{.495\linewidth}
\begin{tikzpicture}
\begin{axis}[width=1.1\linewidth, xlabel={number of terms}, ylabel={accuracy}, xmin=0, xmax=200, ytick={0.5,0.6,0.7,0.8,0.9,1}, yticklabels={0.5,0.6,0.7,0.8,0.9,1}, grid=both, ymin=0.5, ymax=1.0, xlabel near ticks, ylabel near ticks, legend style={draw=none}, legend pos=south east, legend cell align={left}, legend style={row sep=-2.5pt}, clip mode=individual]
\addplot+[solid, thick, Black!80, mark=none] table[x=cost1, y=POD, col sep=comma, mark=none] {./data/kdv_cost.csv};
\addlegendentry{POD}
\addplot+[densely dashed, thick, Blue!80, mark=none] table[x=cost2, y=MAM2, col sep=comma, mark=none] {./data/kdv_cost.csv};
\addlegendentry{$p=2$}
\addplot+[densely dashed, thick, ForestGreen!80, mark=none] table[x=cost3, y=MAM3, col sep=comma, mark=none] {./data/kdv_cost.csv};
\addlegendentry{$p=3$}
\addplot+[densely dashed, thick, Red!80, mark=none] table[x=cost4, y=MAM4, col sep=comma, mark=none] {./data/kdv_cost.csv};
\addlegendentry{$p=4$}
\node[Black!80, mark size=2.5pt] (A) at (14,0.650245) {\pgfuseplotmark{triangle*}};
\node[Blue!80, mark size=2.5pt] (B) at (40,0.860967)  {\pgfuseplotmark{triangle*}};
\node[ForestGreen!80, mark size=2.5pt] (C) at (78,0.935302) {\pgfuseplotmark{triangle*}};
\node[Red!80, mark size=2.5pt] (D) at (128,0.990321) {\pgfuseplotmark{triangle*}};
\end{axis}
\end{tikzpicture}
\caption{}
\end{subfigure}
\caption{Cost-accuracy assessment in evaluating reduced-order model \eqref{eq:rom_nonlinear} for the Korteweg-de Vries problem using (a) the POD-based manifold formulation (Algorithm~\ref{alg:pod-based}) and (b) the alternating minimization based manifold formulation (Algorithm~\ref{alg:am-based}) as a function of the degree of the polynomial embeddings, $p$. The results for the different models at a reduced dimensionality of $r=4$ are indicated by the triangles ($\triangle$).}
\label{fig:kdv_cost}
\end{figure}

When we repeat the experiment but with the dimension of the reduced-order model increased to $r=16$ (choosing $r+q=25$), the state error for the OpInf, MPOD-Opinf, and MAM-OpInf models in the range of the training data drops further to 5.7\%, 4.3\%, and 2.4\% (see Fig.~\ref{fig:kdv_r16}).
However, if the models are integrated to final time $T$, the error increases rapidly.
The state errors across the time interval $[0.2,1]$ are 30.4\%, 48.2\%, and 76.4\%, respectively, for the three OpInf variants.
Fig.~\ref{fig:kdv_snapshots3} shows that these errors can be attributed to inaccurate predictions of the soliton's spatial location. 
These results demonstrate an important tradeoff between choosing $r$, the dimension of the linear subspace described by basis $\mathbf{V}$, and $q$, the number of basis vectors in the basis $\overline{\mathbf{V}}$. 
With $r=16$, the linear subspace captures 99.7\% of the snapshot energy. Enriching the approximation by adding the $\overline{\mathbf{V}}$ manifold terms using the next $q=9$ singular vectors leads to 99.85\% snapshot energy being captured by the polynomial manifold representation (i.e., an increase of only 0.15\%). 
In this case, the components of the basis $\overline{\mathbf{V}}$ correspond to singular vectors with near-zero singular values and the additional terms provide little benefit---in fact, in this example they lead to overfitting and a decline in reduced model predictive performance.

We now shift attention to the computational cost of integrating reduced models of the form \eqref{eq:rom_nonlinear}. 
We monitor the accuracy of the model in the training regime, as given by the representation error of the training data \eqref{eq:rit}, as a function of the length of the reduced-representation vector $[\hat{\mathbf{s}}^\top, (\hat{\mathbf{s}}\otimes\hat{\mathbf{s}})^\top, \hat{\mathbf{g}}(\hat{\mathbf{s}})^\top]^\top$, which equals $r+r(r+1)/2+d(r,p)$. 
The results are summarized in Fig.~\ref{fig:kdv_cost} for the POD-based and AM-based manifold formulations. 
A higher dimension for the reduced-order model leads to both more expensive computation and increased accuracy. 
However, for models of the same dimensionality (for example, $r=4$), it can be seen that the total number of terms grows rapidly with the degree of the polynomial embeddings. 
Note that this analysis pertains only to {\em online} computational costs. 
(Reduced-order modeling of large-scale dynamical systems typically invokes a high cost during the offline phase, which is performed only once.)

\subsection{Incompressible Navier-Stokes -- Flow past a cylinder}
\label{subsec:vortex_shedding}

We now apply our techniques to the well-investigated problem of two-dimensional transient flow past a circular cylinder. 
We focus on the configuration with Reynolds number $Re = 100$, a value that is above the critical Reynolds number for the onset of the two-dimensional vortex shedding. 
The fluid flow is governed by the incompressible Navier-Stokes equations
\begin{equation}
\begin{aligned} 
    \partial_t \mathbf{u} + \nabla \cdot ( \mathbf{u} \otimes \mathbf{u}) &= \nabla p + Re^{-1}\Delta  \mathbf{u}, \\
    \nabla \cdot \mathbf{u} &= 0.
    \label{eq:incom-navier-stokes}
\end{aligned}
\end{equation}
The velocity vector is given by $\mathbf{u}= (u, v)^\top$ where $u$ and $v$ are the components in the $x$ and $y$-direction, respectively. 
Pressure is denoted by $p$. 
We integrate the model over time interval $t \in [0,8]$. 
Problem setup, geometry, and parameters are taken from the DFG 2D-3 benchmark in the FeatFlow benchmark suite.\footnote{\url{https://www.mathematik.tu-dortmund.de/~featflow/en/benchmarks/ff_benchmarks.html}} 
In the model reduction experiments that follow, we did not explicitly account for the pressure term. 
This omission is known to be valid for both the transient and periodic regime of the flow.\cite{deane1991low, noack2003hierarchy}

\begin{figure}[!tbp]
    \centering \footnotesize
\begin{subfigure}{.85\linewidth}
    \begin{overpic}[width=\linewidth]{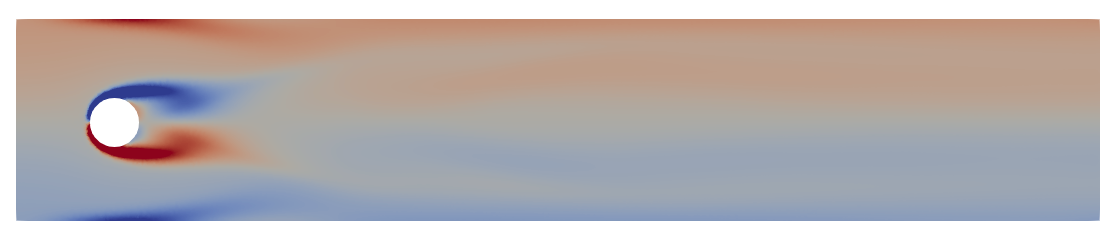}
    \put(75,5){Mean flow}
    \end{overpic}
\end{subfigure}    
\begin{subfigure}{.85\linewidth}
    \begin{overpic}[width=\linewidth]{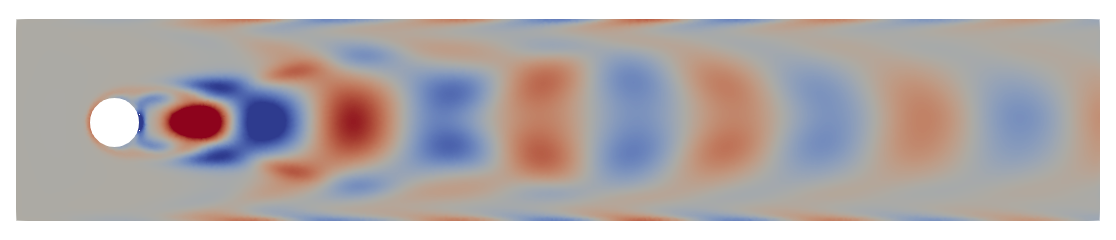}
    \put(75,5){POD mode 1}
    \end{overpic}
\end{subfigure}    
\begin{subfigure}{.85\linewidth}
    \begin{overpic}[width=\linewidth]{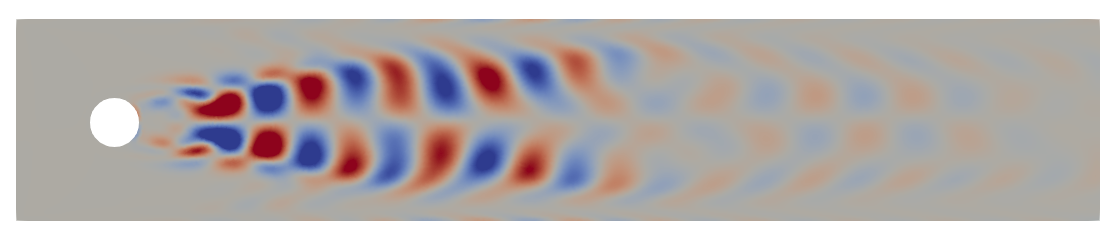}
    \put(75,5){POD mode 3}
    \end{overpic}
\end{subfigure}    
\begin{subfigure}{.85\linewidth}
    \begin{overpic}[width=\linewidth]{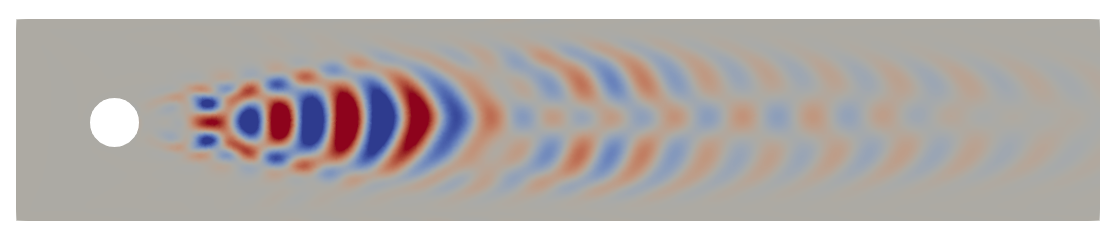}
    \put(75,5){POD mode 5}
    \end{overpic}
\end{subfigure}    
\begin{subfigure}{.85\linewidth}
    \begin{overpic}[width=\linewidth]{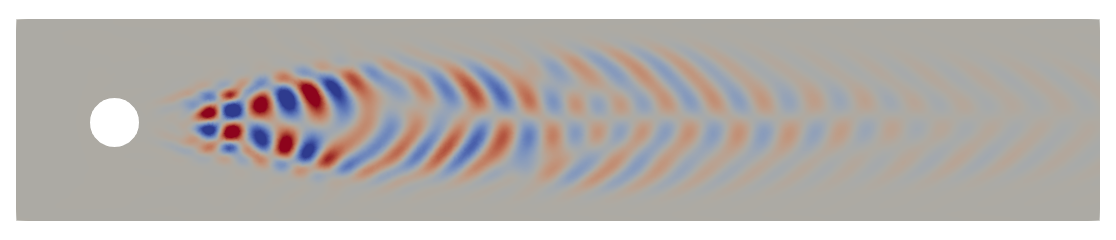}
    \put(75,5){POD mode 7}
    \end{overpic}
\end{subfigure}    
\caption{Mean flow (top) and the dominant odd-numbered POD modes in the cylinder flow problem. We show the vorticity computed from the velocity fields.}
\label{fig:ns_pod_modes}
\end{figure}

We collected 200 snapshots of a periodic reference simulation at $Re=100$ in the interval $t \in [4,5]$, and store each snapshot as a  column vector with 292,678 entries. 
As usual, the snapshot matrix $\mathbf{S}$ is centered by its column-averaged mean value $\mathbf{S}_\text{ref}$, and the orthogonal basis vectors are computed by means of the POD. 
Because the cylinder flow example is periodic, the POD modes can be grouped in pairs $(\mathbf{v}_1, \mathbf{v}_2)$, $(\mathbf{v}_3, \mathbf{v}_4)$, $(\mathbf{v}_5, \mathbf{v}_6)$, $(\mathbf{v}_7, \mathbf{v}_8)$. 
Fig.~\ref{fig:ns_pod_modes} displays the computed mean flow and the first POD mode for each pair.

We choose the reference state in \eqref{eq:nonlinear_approx} to represent the mean flow. 
The flow dynamics can be captured well with only eight modes capturing 99.89\% of the snapshot energy. 
However, physical and mathematical system reduction approaches have revealed that only two modes are {\em actual} degrees of freedom of the system; the remaining ones are completely dependent on these two.\cite{noack2003hierarchy} 
This insight can also be obtained from a nonlinear correlation analysis.\cite{LoiseauBruntonNoack} 
Although the POD analysis indicates that eight POD modes should be considered for accurate flow reconstructions, we use instead  the proposed nonlinear model reduction framework for learning dynamical-system models that respect the problem's intrinsic dimensionality of 2. 
Although we could also compute an orthogonal set of basis vectors through an orthogonal Procrustes problem, as in the alternating minimization based representation learning problem (see Algorithm~\ref{alg:am-based}), the advantages accruing from orthogonal transformations of a POD subspace were found to be negligible: The Procrustes modes were found to be virtually indistinguishable from the ones computed using POD. 
We thus focus exclusively on the MPOD-OpInf formulation.

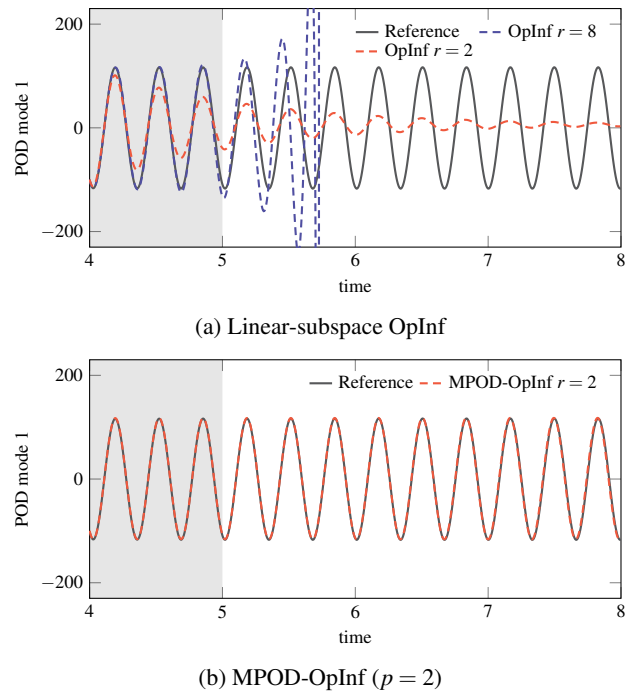
\begin{figure}[!tbp]
\centering \scriptsize
\begin{subfigure}[t]{\linewidth}
\begin{tikzpicture}
\begin{axis}[width=\linewidth, height=.55\linewidth, xlabel={time}, ylabel={POD mode 1}, xmin=4, xmax=8, ymin=-230, ymax=230, grid=none, xlabel near ticks, ylabel near ticks, legend style={draw=none, fill=none}, legend cell align={left}, legend style={row sep=-2.5pt}, legend columns=2]
\fill[black, opacity=0.1] (0.5,-250)node [below=5pt, black] {}--(5,-250)node [below=5pt, black] {}--(5,250)--(3,250)--(3,0)--cycle ;
\addplot+[thick, Black!80, mark=none] table[x=time, y=Shat1x100, col sep=comma, mark=none] {./data/Shat_updated2.csv};
\addlegendentry{Reference}
\addplot+[densely dashed, thick, Blue!80, mark=none] table[x=time, y=shat1, col sep=comma, mark=none] {./data/s_hat_POD_r8.csv};
\addlegendentry{OpInf $r=8$}
\addplot+[densely dashed, thick, Red!80, mark=none] table[x=time, y=shat1, col sep=comma, mark=none] {./data/s_hat_POD_r2.csv};
\addlegendentry{OpInf $r=2$}
\end{axis}
\end{tikzpicture}
\caption{Linear-subspace OpInf}
\label{fig:coef_1}
\end{subfigure}\\[0.5em]
\begin{subfigure}[t]{\linewidth}
\begin{tikzpicture}
\begin{axis}[width=\linewidth, height=.55\linewidth, xlabel={time}, ylabel={POD mode 1}, xmin=4, xmax=8, grid=none, ymin=-230, ymax=230, xlabel near ticks, ylabel near ticks, legend style={draw=none, fill=none}, legend cell align={left}, legend style={row sep=-2.5pt}, legend columns=-1]
\fill[black, opacity=0.1] (0.5,-250)node [below=5pt, black] {}--(5,-250)node [below=5pt, black] {}--(5,250)--(3,250)--(3,0)--cycle ;
\addplot+[thick, Black!80, mark=none] table[x=time, y=Shat1x100, col sep=comma, mark=none] {./data/Shat_updated2.csv};
\addlegendentry{Reference}
\addplot+[densely dashed, thick, Red!80, mark=none] table[x=time, y=shat1, col sep=comma, mark=none] {./data/s_hat_MPOD_r2_p2.csv};
\addlegendentry{MPOD-OpInf $r=2$}
\end{axis}
\end{tikzpicture}
\caption{MPOD-OpInf ($p=2$)}
\label{fig:coef_2}
\end{subfigure}
\caption{Comparison of the amplitudes of the first POD mode over time in the OpInf (top) and the MPOD-OpInf (bottom) models in the cylinder flow problem. The manifold models used quadratic embeddings. The gray shaded region highlights the window over which flow snapshots have been collected for training.}
\label{fig:ns_coefficients}
\end{figure}

\begin{figure}[!tbp]
\centering \scriptsize
\begin{tikzpicture}
\begin{axis}[width=.7\linewidth, height=.7\linewidth, xlabel={POD mode 1}, ylabel={POD mode 2}, xmin=-150, xmax=150, grid=none, ymin=-150, ymax=150, xlabel near ticks, ylabel near ticks, legend style={draw=none}, legend cell align={left}, legend style={row sep=-1.5pt}, legend pos=outer north east, set layers, mark layer=axis tick labels]
\addplot+[color=Black!80,only marks,mark=triangle*,mark options={solid}, each nth point={2}] table[x=Shat1_1cycle, y=Shat2, col sep=comma, Black!80] {./data/Shat.csv};
\addlegendentry{Reference}
\addplot+[densely dashed, thick, Blue!80, mark=none] table[x=shat1_t, y=shat2_t, col sep=comma, mark=none] {./data/s_hat_POD_r8.csv};
\addlegendentry{OpInf $r=8$}
\addplot+[densely dashed, thick, Red!80, mark=none] table[x=shat1, y=shat2, col sep=comma, mark=none] {./data/s_hat_MPOD_r2_p2.csv};
\addlegendentry{MPOD-OpInf $r=2$}
\end{axis}
\end{tikzpicture}
\caption{Limit cycles observed in the full simulation model ($\triangle$), the eight-equation OpInf model (blue curve), and the two-equation MPOD-OpInf model (red curve) for the cylinder flow problem.}
\label{fig:a1-a2}
\end{figure}
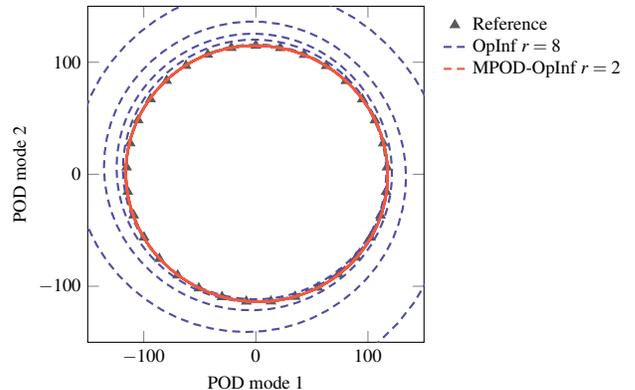

\begin{figure}[!tbp]
\centering \footnotesize
\begin{subfigure}{.9\linewidth}
    \begin{overpic}[width=\linewidth]{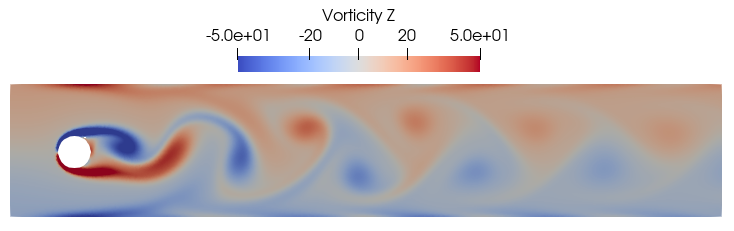}
    \end{overpic}
    \caption{Reference (292,678 degrees of freedom)}
\end{subfigure}    
\begin{subfigure}{.9\linewidth}
    \begin{overpic}[width=\linewidth]{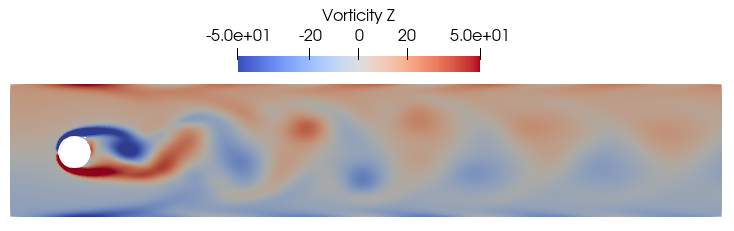}
    \end{overpic}
    \caption{MPOD-OpInf (2 degrees of freedom)}
\end{subfigure}    
\caption{Flow field (vorticity) predicted at $t=8$ with a two-equation MPOD-OpInf model using quadratic embeddings (that is $r=2$; $p=2$) for the cylinder flow problem.}
\label{fig:flow_prediction}
\end{figure}

Fig.~\ref{fig:coef_1} shows that an OpInf model of size $r=2$ is unable to capture the periodic nature of the fluid flow. 
The state error for this model across the training window is 42.9\%, while with all eight modes, the error drops to 7.4\%. 
While the model does an excellent job of capturing the transient dynamics in the training regime $t\in[4,5]$, it fails soon after exiting the training window. 
For the MPOD-OpInf model, we consider only a reduced basis dimension of $r=2$, as informed by physical intuition. 
This means that the remaining six POD modes are contained in the basis matrix associated with the nonlinear part of the state approximation \eqref{eq:nonlinear_approx}, that is, $\overline{\mathbf{V}}$. 
These approximations are built from quadratic embeddings. 
(Although high-order embeddings were considered for this problem, a polynomial degree of $p=2$ was found to be sufficient for learning accurate reduced-order models.) 
The regularization parameter for the representation learning problem was set to $\gamma=10^{-2}$. 
The training error for the MPOD-OpInf was found to be 18.6\%, which is, as expected, larger than the eight-equation OpInf model. 
However, the inferred two-equation MPOD-OpInf model, which has a quadratic term in the nonlinear part of the state approximation, was found to be stable well outside of the training regime (see Fig.~\ref{fig:coef_2}). 
The modal amplitudes of the original simulation model are found when the flow data is projected onto the eight POD modes. 
A comparison in the phase space of the first two coefficients, shown in Fig.~\ref{fig:a1-a2}, finds the MPOD-OpInf model to be accurate and stable with respect to the flow data. 
Finally, Fig.~\ref{fig:flow_prediction} shows a reconstruction of the flow field in the original state space as predicted at time $t=8$ compared to the reference solution at the same time step. 
While some of the finer-scale flow features are not resolved fully, the overall flow dynamics are predicted accurately. 

It should be noted that reduced-order models with excellent predictive performance can be obtained by other means. For instance, approaches in which the models are equipped with \textit{linear} structure (such as the dynamic mode decomposition\cite{schmid2010dynamic, tu2013dynamic, kutz2016dynamic, doi:10.1146/annurev-fluid-030121-015835}) are reported to work well for capturing periodic vortex shedding in the cylinder flow problem.\cite{baddoo2023physics} In the work from Baddoo et al.\cite{baddoo2023physics}, for instance, the data was truncated to the first 15 POD modes. It remains unclear to what extent dynamic mode decomposition and its variants can issue efficient and accurate predictions at or near the \emph{true} dimensionality of two.

\section{Conclusion \& Discussion}
\label{sec:conclusions}

We have presented a general framework for nonlinear model reduction of large-scale physical systems. 
We draw on recently developed techniques for constructing nonlinear manifolds of polynomial structure via representation learning in the form of two different learning approaches.
First, a POD-based version of the approach is intuitive due its connection to conventional POD. 
Second, if one is willing to depart from the interpretable nature of POD methods, alternating minimization techniques can boost model accuracy by means of   better approximations to the solution of the general representation learning problem. 
We then turn to the issue of learning reduced-order models from data. 
By projecting PDE systems onto the nonlinear manifold we can identify the algebraic structure of the projection-based reduced-order model. 
This process calls for careful consideration of the structure of (1) the high-dimensional, physical system and (2) the nonlinear state approximation of choice. 
The non-intrusive OpInf method was used for learning models directly from time-domain simulation data. 
Coupling of the two different representation learning approaches with the OpInf framework leads to a set of methods referred to as POD-based manifold OpInf (MPOD-OpInf) and alternating-minimization-based manifold OpInf (MAM-OpInf).

We applied this methodology to the Allen-Cahn equation, the Korteweg-de Vries equation, and the incompressible Navier-Stokes equation. 
In all numerical experiments, we found the proposed OpInf approaches to be able to circumvent the limitations of linear-subspace OpInf that are due to its use of linear state approximations. 
The polynomial manifold constructions provide the most benefit in situations where the linear subspace does not accurately represent the full dynamics of the training data.
In these situations, the manifold acts as a closure term that accounts for the effects of modes truncated from the linear subspace.
The increased accuracy enabled by a nonlinear compression of the state data does not point to computational speedups: The reduced dimensionality comes at the cost of increased algebraic complexity (and thus computational burden) for the manifold-based reduced models. Although the results from Section IV(D) imply that model robustness and predictive performance are important additional considerations to be made, further investigation is needed in better understanding the tradeoffs between dimensionality and complexity.

Further improvements in OpInf reduced-order models may be possible if constraints are introduced to enforce particular mathematical properties of the dynamical system. 
For example, some classes of problems can be expressed using Hamiltonian or Lagrangian formalisms.\cite{SHARMA2022133122} 
Biasing OpInf models toward such structure may enable more accurate long-time predictions far outside the training time interval and will be addressed in future work. In another research direction, data-driven OpInf could be combined with the dynamic training via roll outs of neural ordinary differential equations.\cite{UY2023224} This OpInf formulation should be more robust against perturbations in the data because the whole predicted trajectory is considered in the training loss rather than a single time step.

\begin{acknowledgments}
This work has been supported in part by the U.S.\ Department of Energy AEOLUS MMICC center under award DE-SC0019303, program manager W.\ Spotz, and by the AFOSR MURI on physics-based machine learning, award FA9550-21-1-0084, program manager F.\ Fahroo. 
Stephen Wright was also supported by NSF awards DMS 2023239 and CCF 2224213.
Laura Balzano was supported by ARO YIP award W911NF1910027 and NSF CAREER award CCF-1845076.
Karen Willcox would also like to thank Yvon Maday and Albert Cohen for several helpful conversations during the 2022 Leçons Jacques-Louis Lions in Paris, France.
\end{acknowledgments}

\section*{Data Availability Statement}

A Jupyter notebook outlining the representation learning problem for inferring, from data, nonlinear state approximations of the form \eqref{eq:nonlinear_approx} for the problem from Section~\ref{subsec:illustrative_example} is available at \url{https://github.com/geelenr/nl_manifolds}. The notebook features both the POD and alternating minimization based formulations from Algorithms \ref{alg:pod-based} and \ref{alg:am-based}. The data used in numerical experiments from Section~\ref{subsec:kdv} and \ref{subsec:allen-cahn} are available upon reasonable request from the authors. The FEniCSx computing platform is used to solve the equations \eqref{eq:incom-navier-stokes} through their tutorial example.\cite{alnaes2015fenics, logg2012automated}

\section*{References}

\bibliography{aipsamp}% Produces the bibliography via BibTeX.
\bibliographystyle{naturemag}

\end{document}